\documentclass[11pt]{article} 
\usepackage{latexsym} 
\newcommand{\eproof}{\mbox{\ }\hfill $\Box$ \par \vskip 10pt}

\newtheorem{Theorem}{Theorem}[section] 
\newtheorem{lemma}[Theorem]{Lemma} 
\newtheorem{prop}[Theorem]{Proposition}

\begin{document}

\title{Dispersive estimates of solutions to the
Schr\"odinger equation in dimensions $n\ge 4$}

\author{{\sc Georgi Vodev}}

\date{} 
 
\maketitle 

\begin{abstract}
\noindent
We prove dispersive estimates 
for solutions to the Schr\"odinger  
equation with a real-valued potential
$V\in L^\infty({\bf R}^n)$, $n\ge4$, satisfying 
$V(x)=O(\langle x\rangle^{-(n+2)/2-\epsilon})$, $\epsilon>0$. 
\end{abstract}

\setcounter{section}{0}
\section{Introduction and statement of results}

Let $V\in L^\infty({\bf R}^n)$, $n\ge 4$,
 be a real-valued function satisfying
$$\left|V(x)\right|\le C\langle x\rangle^{-\delta},
\quad \forall x\in {\bf R}^n,\eqno{(1.1)}$$
with constants $C>0$ and $\delta>(n+2)/2$, 
where $\langle x\rangle=(1+|x|^2)^{1/2}$. 
Denote by $G_0$ and $G$ the self-adjoint realizations of the
operators $-\Delta$ and $-\Delta+V(x)$ on $L^2({\bf R}^n)$.
It is well known that the absolutely continuous spectrums of the
operators $G_0$ and $G$ coincide with the interval $[0,+\infty)$.
Moreover, by Kato's theorem the operator $G$ has no strictly 
positive eigenvalues, which in turn implies that  
$G$ has no strictly positive resonances neither. 
Throughout this paper, given $1\le p\le +\infty$, $L^p$ will 
denote the space $L^p({\bf R}^n)$. Also, 
given an $a>0$ denote by $\chi_a\in C^\infty({\bf R})$ 
 a real-valued function supported in the interval
$[a,+\infty)$, $\chi_a=1$ on $[2a,+\infty)$. 
Our main result is the following

\begin{Theorem}
Assume (1.1) fulfilled. Then 
for every $a>0$, $0<\epsilon\ll 1$,   
there exist constants $C, C_\epsilon>0$ so that the following estimates hold
$$\left\|e^{itG}G^{-(n-3)/4}
\chi_a(G)\right\|_{L^1\to L^\infty}
\le C|t|^{-n/2},\quad \forall t\neq 0,\eqno{(1.2)}$$
$$\left\|e^{itG}
\chi_a(G)\langle x\rangle^{-n/2-\epsilon}
\right\|_{L^2\to L^\infty}\le C_\epsilon|t|^{-n/2},
\quad \forall t\neq 0.\eqno{(1.3)}$$
Moreover, for every 
 $0\le q\le (n-3)/2$, $2\le p<\frac{2(n-1-2q)}{n-3-2q}$, we have
$$\left\|e^{itG}G^{-\alpha q/2}
\chi_a(G)\right\|_{L^{p'}\to L^p}
\le C|t|^{-\alpha n/2},\quad \forall t\neq 0,\eqno{(1.4)}$$
where $1/p+1/p'=1$, $\alpha=1-2/p$. 
\end{Theorem}

{\bf Remark 1.} The desired result would be to prove the estimate
$$\left\|e^{itG}\chi_a(G)\right\|_{L^1\to L^\infty}
\le C|t|^{-n/2},\quad \forall t\neq 0.\eqno{(1.5)}$$
It is shown by Goldberg and Visan \cite{kn:GV}, however, that when
$n\ge 4$ there exists a compactly supported potential
$V\in C^k({\bf R}^n)$, $\forall k<(n-3)/2$, for which (1.5) fails to hold. 
In other words, in order that (1.5) holds true one needs to have a control
 of $(n-3)/2$ derivatives of $V$. It seems that
for potentials satisfying (1.1) only, our estimate (1.2) with a loss
of $(n-3)/2$ derivatives is quite optimal.

{\bf Remark 2.} It is natural to expect that if the zero is
neither an eigenvalue nor a resonance of $G$, the statements
of Theorem 1.1 hold true with $\chi_a$ replaced by the
characteristic function, $\chi$, of the interval $[0,+\infty)$
(the absolutely continuous spectrum of $G$),  
$G^{-(n-3)/4}$ and $G^{-\alpha q/2}$ replaced by $\left\langle G
\right\rangle^{-(n-3)/4}$ and $\left\langle G
\right\rangle^{-\alpha q/2}$, respectively. 
To prove this, it suffices to show
that in this case the estimate (1.5) holds true with $\chi_a$
replaced by $(1-\chi_a)\chi$ for $a>0$ small enough. The proof of
such an estimate, however, requires different techniques than those
developed in the present paper. 

{\bf Remark 3.} We believe that the estimate (1.4) still holds at
the end point $p=\frac{2(n-1-2q)}{n-3-2q}$ for $0\le q<(n-3)/2$
and that our approach leads to such an estimate. 
In fact, it is not hard to see from
the proof of (1.4) in Section 4 that the problem is reduced to
estimating the $L^2\to L^2$ norm of operators with explicitly
given kernels.

The estimate (1.5) is proved in the case $n=2$ by Schlag \cite{kn:S}
for potentials satisfying (1.1) with $\delta>2$. When $n=3$, 
the estimate (1.5) is proved by Goldberg and Schlag \cite{kn:GS} 
for potentials satisfying (1.1) with $\delta>3$. Recently, 
(1.5) has been proved in this case 
by Vodev \cite{kn:V1} and Yajima \cite{kn:Y2}
for potentials satisfying (1.1) with $\delta>5/2$, while 
Goldberg \cite{kn:G} has proved (1.5) for a very large class of
potentials including those satisfying (1.1) with $\delta>2$.
The proofs in all these papers (except for \cite{kn:V1}) are based
on the very nice properties of the outgoing and incoming 
free resolvents when $n=2$ and $n=3$. 
When $n\ge 4$, however, these properties are no longer valid, 
and consequently one needs different methods to prove estimates like
(1.5) (or like these in Theorem 1.1). The first result in this case
is due to Journ\'e, Sofer and Sogge \cite{kn:JSS}, where
they proved an analogue of (1.5) for potentials satisfying (1.1)
with $\delta>n+4$ as well as the regularity property 
$\hat V\in L^1$. This was later improved by Yajima \cite{kn:Y1}
using the properties of the wave operators. Note also that
the estimate (1.3) was proved by Jensen and Nakamura \cite{kn:JN}
for potentials satisfying (1.1) with $\delta>n$ as well as an
extra technical assumption.

To prove Theorem 1.1 we use the method of \cite{kn:V1} together with
some ideas from \cite{kn:V2} where similar dispersive estimates
have been proved for the wave group $e^{it\sqrt{G}}$ for 
potentials satisfying (1.1) with $\delta>(n+1)/2$. Roughly
speeking, the method consists of reducing the dispersive estimates
to uniform estimates for the Shr\"odinger group (resp. the wave group)
on weihted $L^2$ spaces, which in turn are proved by using some
more or less known properties of the perturbed resolvent on 
 weihted $L^2$ spaces (see Section 3). Note finally that in view of
Goldberg's result \cite{kn:G}, one should expect that the
statements of Theorem 1.1 hold true for the larger class of
potentials satisfying (1.1) with $\delta>(n+1)/2$. The proof of
such estimates, however, would require a different approach than
this one presented here.

\section{Preliminary estimates}

The following properties of the free Schr\"odinger group 
will play a key role in the proof of our dispersive estimates.

\begin{prop}
Let $\psi\in C_0^\infty((0,+\infty))$. For every 
$0\le s\le (n-1)/2$, $0<\epsilon\ll 1$, $0<h\le 1$, 
$t\neq 0$, we have 
$$\left\|e^{itG_0}\psi(h^2G_0)
\langle x\rangle^{-1/2-s-\epsilon}\right\|_{L^2\to L^\infty}\le
 Ch^{-(n-1)/2+s}|t|^{-s-1/2},\eqno{(2.1)}$$
with a constant $C>0$ independent of $t$ and $h$.
For every $s\ge 0$, $0<h\le 1$, $t\in {\bf R}$, we have
$$\left\|\langle x\rangle^{-s}
e^{itG_0}\psi(h^2G_0)\langle x\rangle^{-s}
\right\|_{L^2\to L^2}\le C\langle t/h\rangle^{-s},\eqno{(2.2)}$$
with a constant $C>0$ independent of $t$ and $h$.
\end{prop}

{\it Proof.} We are going to take advantage of the formula
$$e^{itG_0}\psi(h^2G_0)=(\pi i)^{-1}\int_0^\infty
 e^{it\lambda^2}\psi(h^2\lambda^2)\left(R_0^+(\lambda)-R_0^-(\lambda)\right)
\lambda d\lambda,\eqno{(2.3)}$$
where $R_0^\pm(\lambda)$ are the outgoing and incoming free resolvents
with kernels given in terms of the Hankel functions by
$$[R^\pm_0(\lambda)](x,y)=\pm\frac{i}{4}\left(\frac{\lambda}
{2\pi |x-y|}\right)^\nu H_\nu^\pm(\lambda|x-y|),$$
where $\nu=(n-2)/2$. Hence, the kernel of the operator
$e^{itG_0}\psi(h^2G_0)$ is of the form $K_h(|x-y|,t)$, where
$$K_h(\sigma,t)=\frac{\sigma^{-2\nu}}{(2\pi)^{\nu+1}}
\int_0^\infty
 e^{it\lambda^2}\psi(h^2\lambda^2){\cal J}_\nu(\sigma\lambda)
\lambda d\lambda=h^{-n}K_1(\sigma h^{-1},th^{-2}),\eqno{(2.4)}$$
where ${\cal J}_\nu(z)=z^\nu J_\nu(z)$, $J_\nu(z)=(H^+_\nu(z)+
H^-_\nu(z))/2$ is the Bessel function of order $\nu$. 
It is easy to see that (2.1) follows from the bound
$$\left|K_h(\sigma,t)\right|\le Ch^{s-(n-1)/2}|t|^{-s-1/2}\sigma^{-(n-1)/2+s},
\quad \forall t\neq 0,\sigma>0, 0<h\le 1.\eqno{(2.5)}$$
In view of (2.4), it suffices to show (2.5) with $h=1$.
Let $m\ge 0$ be an integer. Integrating $m$ times in (2.4)
we obtain
$$(it)^mK_1(\sigma,t)=\frac{\sigma^{-2\nu}}{2(2\pi)^{\nu+1}}
\int_0^\infty
 e^{it\mu}\frac{d^m}{d\mu^m}\left(\psi(\mu){\cal J}_\nu(\sigma\sqrt{\mu})
\right) d\mu$$ 
 $$=\sum_{k=0}^m\sigma^{k-2\nu}\int_0^\infty
 e^{it\mu}\psi_{m-k}(\mu){\cal J}_\nu^{(k)}(\sigma\sqrt{\mu})d\mu,$$
where $\psi_{m-k}\in C_0^\infty((0,+\infty))$, 
${\cal J}_\nu^{(k)}(z):=d^k{\cal J}_\nu(z)/dz^k.$
Making the change $\mu=\lambda^2$ we can write the above identity 
in the form
$$(it)^mK_1(\sigma,t)=\sum_{k=0}^m\sigma^{k-2\nu}\int_0^\infty
 e^{it\lambda^2}\varphi_{m-k}(\lambda){\cal J}_\nu^{(k)}(\sigma\lambda)
d\lambda,\eqno{(2.6)}$$
where $\varphi_{k}(\lambda)=2\lambda\psi_k(\lambda^2)$. To estimate
the integral in the RHS of (2.6), we will use that, given any $b>0$ 
and a function $\varphi\in C_0^\infty([-b,b])$, we have the bound
(see the proof of Lemma 2.4 of \cite{kn:RS}):
$$\left|\int_{-\infty}^\infty e^{it\lambda^2-i\tau\lambda}\varphi(
\lambda)d\lambda\right|\le C|t|^{-1/2}\|\hat\varphi\|_{L^1}$$
 $$\le C_b|t|^{-1/2}\sup_{0\le j\le 1}\sup_{\lambda\in
{\bf R}}\left|\partial_\lambda^j\varphi(\lambda)\right|,
\quad\forall t\neq 0,\,\tau\in {\bf R},\eqno{(2.7)}$$
with a constant $C_b>0$ independent of $t$, $\tau$ and $\varphi$.
Consider first the case $0<\sigma\le 1$. 
It is well known that near $z=0$ the function ${\cal J}_\nu(z)$ 
 is equal to $z^{2\nu}$ times an analytic function.
Therefore, we have, for $0<z\le z_0$,
 $$\left|\partial_z^k{\cal J}_\nu(z)\right|\le 
Cz^{2\nu-k},\eqno{(2.8)}$$
for all integers $0\le k\le n-2$, while for integers $k\ge n-1$,
 $$\left|\partial_z^k{\cal J}_\nu(z)\right|\le C_k,\eqno{(2.9)}$$
with constants $C,C_k>0$ depending on $z_0$. 
By (2.6)-(2.9), we obtain
$$\left|K_1(\sigma,t)\right|\le C_m|t|^{-m-1/2},\quad\forall t\neq 0,\,
0<\sigma\le 1,\eqno{(2.10)}$$
for all integers $m\ge 0$ with a constant $C_m>0$ independent
of $t$ and $\sigma$. Clearly, (2.10) holds with $m=s$ for all
real $0\le s\le (n-1)/2$, which in turn implies (2.5) with $h=1$
in this case.

Let now $\sigma\ge 1$.
It is well known that ${\cal J}_\nu(z)=e^{iz}b^+_\nu(z)+
e^{-iz}b^-_\nu(z)$ with functions $b^\pm_\nu(z)$ satisfying
$$\left|\partial_z^jb^\pm_\nu(z)\right|\le C_jz^{(n-3)/2-j},
\quad\forall z\ge z_0,\eqno{(2.11)}$$
for every integer $j\ge 0$ and every $z_0>0$, with a constant
$C_j>0$ depending on $j$ and $z_0$ but independent of $z$.
We can write $K_1=K^+_1+K^-_1$, where $K_1^\pm$ is defined by replacing in
(2.4) the function ${\cal J}_\nu(z)$ by $e^{\pm iz}b^\pm_\nu(z)$.
We have
$${\cal J}_\nu^{(k)}(z)=\sum_\pm e^{\pm iz}b^\pm_{\nu,k}(z)\eqno{(2.12)}$$
with functions $b^\pm_{\nu,k}$ satisfying (2.11). Thus, by 
(2.6), (2.7), (2.11) and (2.12) we get
$$\left|K_1^\pm(t,\sigma)\right|\le C_m|t|^{-m-1/2}\sigma^{-(n-1)/2+m},
\quad\forall t\neq 0,\,\sigma\ge 1,\eqno{(2.13)}$$
for all integers $m\ge 0$ with a constant $C_m>0$ independent
of $t$ and $\sigma$. Obviously, (2.13) holds true with $m=s$ for all 
real $0\le s\le (n-1)/2$, which in turn implies (2.5) with $h=1$
in this case.

Given a set ${\cal M}\subset{\bf R}^n$ denote by
$\eta({\cal M})$ the characteristic function of ${\cal M}$. We have
$$\left\|\langle x\rangle^{-s}
e^{itG_0}\psi(h^2G_0)\langle x\rangle^{-s}
\right\|_{L^2\to L^2}$$ $$\le \left\|\eta(|x|\le \gamma|t|/2h)
e^{itG_0}\psi(h^2G_0)\eta(|x|\le \gamma|t|/2h)
\right\|_{L^2\to L^2}+C\langle t/h\rangle^{-s},\eqno{(2.14)}$$
where $\gamma>0$ is a constant to be fixed below. 
In view of Schur's lemma the norm in the RHS of (2.14) is upper bounded
by
$$\sup_{|x|\le\gamma|t|/2h}\int_{|y|\le\gamma|t|/2h}|K_h(|x-y|,t)|dy\le
\int_{|\xi|\le\gamma|t|/h}|K_h(|\xi|,t)|d\xi$$ $$\le
 C\int_0^{\gamma|t|/h}\sigma^{n-1}\left|K_h(\sigma,t)
\right|d\sigma.\eqno{(2.15)}$$
Thus, to prove (2.2) it suffices to show that the integral in the RHS
of (2.15) is upper bounded by $O_m\left(\langle t/h\rangle^{-m}\right)$
for all integers $m\ge 0$. In view of (2.4) this would follow from
the bound
$$\int_0^{\gamma|t|}\sigma^{n-1}\left|K_1(\sigma,t)
\right|d\sigma\le C_m\langle t\rangle^{-m},\eqno{(2.16)}$$
provided $\gamma$ is properly choosen. Write the function
$K_1^\pm$ in the form
$$K_1^\pm(\sigma,t)=\frac{\sigma^{-2\nu}e^{-i\sigma^2/4t}}{(2\pi)^{\nu+1}}
\int_0^\infty
 e^{it(\lambda\pm\sigma/2t)^2}\varphi(\lambda)b^\pm_\nu(\sigma\lambda)
 d\lambda,\eqno{(2.17)}$$
where $\varphi(\lambda)=\lambda\psi(\lambda^2)$. Clearly, we can fix now
$\gamma>0$ (depending on supp$\,\varphi$) so that
$|\lambda\pm\sigma/2t|\ge Const>0$ on supp$\,\varphi$ for 
$\sigma\le\gamma|t|$. Therefore, integrating by parts $m$ times in
(2.17) and using (2.11), one can easily obtain the bound
$$\left|K_1^\pm(\sigma,t)\right|\le C_m|t|^{-m},\quad 1\le \sigma\le
\gamma|t|,\eqno{(2.18)}$$
for every integer $m\ge 0$. On the other hand, in the same way as in
the proof of (2.10), one concludes that $K_1$ satisfies (2.18)
for $0<\sigma\le 1$ as well. This clearly implies (2.16), and hence
(2.2).
\eproof

We will also need the following lemma proved in \cite{kn:V2} 
(see also \cite{kn:V1}).

\begin{lemma} Assume (1.1) fulfilled. Then, for every 
$\psi\in C_0^\infty((0,+\infty))$, $0\le s\le \delta$, 
$1\le p\le \infty$, $0<h\le 1$, we have
$$\left\|\langle x\rangle^{-s}\psi(h^2G_0)
\langle x\rangle^{s}\right\|_{L^2\to L^2}\le C,\eqno{(2.19)}$$
 $$\left\|\langle x\rangle^{-s}\psi(h^2G)
\langle x\rangle^{s}\right\|_{L^2\to L^2}\le C,\eqno{(2.20)}$$
 $$\left\|\left(\psi(h^2G_0)
-\psi(h^2G)\right)
\langle x\rangle^{s}\right\|_{L^2\to L^2}\le Ch^2,\eqno{(2.21)}$$
 $$\left\|\psi(h^2G_0)\right\|_{L^p\to L^p}\le C,\eqno{(2.22)}$$
 $$\left\|\psi(h^2G)\right\|_{L^p\to L^p}\le C,\eqno{(2.23)}$$
 $$\left\|\psi(h^2G)-\psi(h^2G_0)
\right\|_{L^p\to L^p}\le Ch^2,\eqno{(2.24)}$$
 $$\left\|\psi(h^2G_0)\right\|_{L^2\to L^p}\le Ch^{-n|\frac{1}{2}
-\frac{1}{p}|},\eqno{(2.25)}$$
 $$\left\|\psi(h^2G)\right\|_{L^2\to L^p}\le C
h^{-n|\frac{1}{2}-\frac{1}{p}|},\eqno{(2.26)}$$
 $$\left\|\left(\psi(h^2G)-\psi(h^2G_0)\right)
\langle x\rangle^{s}\right\|_{L^2\to L^p}\le C
h^{2-n|\frac{1}{2}-\frac{1}{p}|},\eqno{(2.27)}$$
with a constant $C>0$ independent of $h$.
\end{lemma}

\section{$L^2\to L^2$ estimates for the Schr\"odinger group}

Given a parameter $0<h\le 1$, and a real-valued 
function $\psi\in C_0^\infty((0,+\infty))$, denote
$$\Psi(t;h)=e^{itG}\psi(h^2G)-
e^{itG_0}\psi(h^2G_0).$$
We will first prove the following

\begin{Theorem} Assume (1.1) fulfilled. Then, we have
$$\left\|\Psi(t;h)\right\|_{L^2\to L^2}\le Ch,\quad
\forall t,\,0<h\le 1,\eqno{(3.1)}$$
with a constant $C>0$ independent of $t$ and $h$.
\end{Theorem}

{\it Proof.} We will derive (3.1) from the following

\begin{prop} For every $s>1/2$, $0<h\le 1$, we have
$$\int_{-\infty}^\infty\left\|\langle x\rangle^{-s}e^{itG}
\psi(h^2G)f\right\|^2_{L^2}dt\le Ch\|f\|^2_{L^2},\quad\forall
 f\in L^2,\eqno{(3.2)}$$
with a constant $C>0$ independent of $f$ and $h$.
\end{prop}

Let $\psi_1\in C_0^\infty((0,+\infty))$ be a real-valued function
 such that $\psi_1\psi\equiv 
\psi$. By Duhamel's formula we obtain the identity
$$\Psi(t;h)=\sum_{j=1}^2\Psi_j(t;h),\eqno{(3.3)}$$
where
$$\Psi_1(t;h)=\left(\psi_1(h^2G)-\psi_1(h^2G_0)\right)\Psi(t;h)$$
$$+\left(\psi_1(h^2G)-\psi_1(h^2G_0)\right)e^{itG_0}\psi(h^2G_0)
-\psi_1(h^2G_0)e^{itG_0}\left(\psi(h^2G)-\psi(h^2G_0)\right),$$
$$\Psi_2(t;h)=i\int_0^t e^{i(t-\tau)G_0}\psi_1(h^2G_0)Ve^{i
\tau G}\psi(h^2G)d\tau.$$
In view of (2.21) we have
$$\left\|\Psi_1(t;h)\right\|_{L^2\to L^2}\le Ch^2.\eqno{(3.4)}$$
For all nontrivial $f,g\in L^2$,  
we have with $0<s-1/2\ll 1$, $\forall\gamma>0$,
$$\left|\left\langle \Psi_2(t;h)f,g\right\rangle\right|\le 
\int_{-\infty}^\infty\left|\left\langle \langle x\rangle^{s}V
 e^{i\tau G}\psi(h^2G)f,\langle x\rangle^{-s} e^{i(t-\tau)G_0}\psi_1(h^2G_0)g
\right\rangle\right|d\tau$$ $$
\le C\gamma\int_{-\infty}^\infty\left\|\langle x\rangle^{-s}
 e^{i\tau G}\psi(h^2G)f\right\|^2_{L^2}d\tau 
 +C\gamma^{-1}\int_{-\infty}^\infty\left\|\langle x\rangle^{-s}
 e^{i\tau G_0}\psi_1(h^2G_0)g
\right\|^2_{L^2}d\tau$$ $$\le  
Ch\gamma\|f\|^2_{L^2}+
 Ch\gamma^{-1}\|g\|^2_{L^2}\le 
 O(h)\|f\|_{L^2}\|g\|_{L^2},\eqno{(3.5)}$$
if we choose $\gamma=\|g\|_{L^2}/\|f\|_{L^2}$, 
which implies (3.1).
\eproof

{\it Proof of Proposition 3.2.} 
Denote by ${\cal H}$ the Hilbert space
$L^2({\bf R};L^2)$. Clearly, (3.2) is equivalent to the fact that
the operator ${\cal A}_h:L^2\to{\cal H}$ defined by
$$\left({\cal A}_hf\right)(x,t)=\langle x\rangle^{-s}e^{itG}
\psi(h^2G)f$$
is bounded with norm $O(h^{1/2})$. Observe that the adjoint
${\cal A}_h^*:{\cal H}\to L^2$ is defined by
$${\cal A}_h^*f=\int_{-\infty}^\infty e^{-i\tau G}\psi(h^2G)
\langle x\rangle^{-s}f(\tau,x)d\tau,$$
so we have, $\forall f,g\in{\cal H}$,
$$\left\langle {\cal A}_h{\cal A}_h^*f,g\right\rangle_{{\cal H}}=
\int_{-\infty}^\infty\left\langle \rho(t,\cdot),g(t,\cdot)
\right\rangle_{L^2}dt,\eqno{(3.6)}$$
where
$$\rho(t,x)=\int_{-\infty}^\infty \langle x\rangle^{-s}
e^{i(t-\tau)G}\psi^2(h^2G)
\langle x\rangle^{-s}f(\tau,\cdot)d\tau.$$
Hence, for the Fourier transform, $\hat\rho(\lambda,x)$, of
$\rho(t,x)$ with respect to the variable $t$ we have
$$\hat\rho(\lambda,x)=Q(\lambda)\hat f(\lambda,x),\eqno{(3.7)}$$
where $Q(\lambda)$ is the Fourier transform of the operator
$$\langle x\rangle^{-s}
e^{itG}\psi^2(h^2G)\langle x\rangle^{-s}.$$
On the other hand, the formula
$$e^{itG}\psi^2(h^2G)=\frac{1}{\pi i}\int_{-\infty}^\infty
e^{it\lambda^2}\psi^2(h^2\lambda^2)\left(R^+(\lambda)-R^-(\lambda)\right)
\lambda d\lambda,\eqno{(3.8)}$$
where
$$R^\pm(\lambda)=\lim_{\varepsilon\to 0^+}(G-\lambda^2\pm i
\varepsilon)^{-1}:\langle x\rangle^{-s}L^2\to 
\langle x\rangle^{s}L^2,\quad s>1/2,$$
shows that
$$Q(\lambda)=(2\pi i)^{-1}\psi^2(h^2\lambda)
\langle x\rangle^{-s}\left(R^+(\sqrt{\lambda})-R^-(\sqrt{\lambda})\right)
\langle x\rangle^{-s}.\eqno{(3.9)}$$
Note that the limit exists in view of the limiting absorption principle.
Moreover, we have the estimate (e.g. see Lemma 3.3 of \cite{kn:V2})
$$\left\|\langle x\rangle^{-s}R^\pm(\lambda)
\langle x\rangle^{-s}\right\|_{L^2\to L^2}\le C\lambda^{-1},
\quad \lambda\ge\lambda_0,\eqno{(3.10)}$$
for every $s>1/2$, $\lambda_0>0$, with a constant
$C>0$ independent of $\lambda$. By (3.9) and (3.10) we conclude
$$\|Q(\lambda)\|_{L^2\to L^2}\le Ch\eqno{(3.11)}$$
with a constant $C>0$ independent of $\lambda$ and $h$. By (3.7) and
(3.11),
$$\|\hat\rho(\lambda,\cdot)\|_{L^2}\le Ch\|\hat f(\lambda,\cdot)\|_{L^2},
\eqno{(3.12)}$$
which together with (3.6) leads to
$$\left|\left\langle {\cal A}_h{\cal A}_h^*f,g\right
\rangle_{{\cal H}}\right|=\left|
\int_{-\infty}^\infty\left\langle \hat\rho
(\lambda,\cdot),\hat g(\lambda,\cdot)
\right\rangle_{L^2}d\lambda\right|$$ $$\le Ch
\int_{-\infty}^\infty\|\hat f
(\lambda,\cdot)\|_{L^2}\|\hat g(\lambda,\cdot)\|_{L^2}d\lambda$$ 
 $$\le Ch\gamma\int_{-\infty}^\infty\|\hat f
(\lambda,\cdot)\|_{L^2}^2d\lambda+
Ch\gamma^{-1}\int_{-\infty}^\infty\|\hat g
(\lambda,\cdot)\|_{L^2}^2d\lambda$$ 
 $$=Ch\gamma\|f\|^2_{{\cal H}}+Ch\gamma^{-1}\|g\|^2_{{\cal H}}=2Ch
\|f\|_{{\cal H}}\|g\|_{{\cal H}},\eqno{(3.13)}$$
if we take $\gamma=\|g\|_{{\cal H}}/\|f\|_{{\cal H}}$, with a constant
$C>0$ independent of $h$, $f$ and $g$. It follows from (3.13) that
the operator ${\cal A}_h{\cal A}_h^*:{\cal H}\to{\cal H}$ is bounded
with norm $O(h)$, and hence the operator ${\cal A}_h:L^2\to
{\cal H}$ is bounded with norm $O(h^{1/2})$.
\eproof

In what follows in this section we will prove the following

\begin{Theorem} Assume (1.1) fulfilled. Then, for every real-valued function 
$\psi\in C_0^\infty((0,+\infty))$ and every 
$0\le s\le n/2$, $0<\epsilon\ll 1$, we have
$$\left\|\langle x\rangle^{-s-\epsilon}
e^{itG}\psi(h^2G)\langle x\rangle^{-s-\epsilon}
\right\|_{L^2\to L^2}\le C\langle t/h\rangle^{-s},\quad 
\forall t,\,0<h\le 1,\eqno{(3.14)}$$
with a constant $C>0$ independent of $t$ and $h$.
\end{Theorem}

{\it Proof.} We will derive (3.14) from the following estimates 

\begin{prop} Assume (1.1) fulfilled. Then, for every real-valued function 
$\psi\in C_0^\infty((0,+\infty))$ and every 
$0\le s\le n/2$, $0<\epsilon\ll 1$, $0<h\le 1$, we have
$$\left\|\langle x\rangle^{-1/2-s-\epsilon}
e^{itG}\psi(G)\langle x\rangle^{-1/2-s-\epsilon}
\right\|_{L^2\to L^2}\le C\langle t\rangle^{-s},\quad 
\forall t,\eqno{(3.15)}$$
$$\int_{-\infty}^\infty\langle t/h\rangle^{2s}
\left\|\langle x\rangle^{-1/2-s-\epsilon}e^{itG}
\psi(h^2G)\langle x\rangle^{-1/2-s-\epsilon}
f\right\|^2_{L^2}dt\le Ch\|f\|^2_{L^2},\quad\forall f\in L^2,\eqno{(3.16)}$$
with a constant $C>0$ independent of $t$, $h$ and $f$. 
\end{prop}

By a standard interpolation argument (e.g. see the proof of Theorem 1.2
of \cite{kn:V1} or the proof of Theorem 3.4 of \cite{kn:V2}) one can 
easily conclude that it suffices to prove (3.14) with 
$s=n/2$, only. On the other hand, in view of (2.2), it suffices to
prove (3.14) for the difference $\Psi(t;h)$. To do so, we will make
use of (3.3). Using (1.1), (2.2) and (2.21), we obtain, 
$$\left\|\langle x\rangle^{-s-\epsilon}\Psi_1(t;h)
\langle x\rangle^{-s_1-\epsilon}f\right\|_{L^2}$$ $$\le O(h^2)
\left\|\langle x\rangle^{-(n+2)/2-\epsilon}\Psi(t;h)
\langle x\rangle^{-s_1-\epsilon}f\right\|_{L^2}+O(h^2)
\langle t/h\rangle^{-s}\|f\|_{L^2},\quad\forall f\in L^2,\eqno{(3.17)}$$
for all $s_1\ge s\ge 0$. 
Using (1.1), (2.2) and (3.16), we obtain, $\forall f,g\in L^2$, with
$s=n/2$,
$$\langle t/h\rangle^s
\left|\left\langle\langle x\rangle^{-s-\epsilon}\Psi_2(t;h)
\langle x\rangle^{-s-1/2-\epsilon}f,g\right\rangle\right|$$ $$\le
\langle t/h\rangle^s
\int_0^t\left|\left\langle\langle x\rangle^{-s-\epsilon}e^{i(t-\tau)G_0}
\psi_1(h^2G_0)Ve^{i\tau G}\psi(h^2G)
\langle x\rangle^{-s-1/2-\epsilon}f,g\right\rangle\right|d\tau$$
 $$\le C\int_0^{t/2}
\langle (t-\tau)/h\rangle^s\left\|\langle x\rangle^{-n/2-
\epsilon}e^{-i(t-\tau)G_0}
\psi_1(h^2G_0)\langle x\rangle^{-s-\epsilon}g\right\|_{L^2}$$
 $$\times \left\|\langle x\rangle^{-1-\epsilon}
e^{i\tau G}\psi(h^2G)\langle x\rangle^{-s-1/2-\epsilon}f
\right\|_{L^2}d\tau$$
  $$+C\int_{t/2}^t\left\|\langle x\rangle^{-1/2-
\epsilon}e^{-i(t-\tau)G_0}
\psi_1(h^2G_0)\langle x\rangle^{-s-\epsilon}g\right\|_{L^2}$$
 $$\times \langle\tau/h\rangle^s
\left\|\langle x\rangle^{-(n+1)/2-\epsilon}
e^{i\tau G}\psi(h^2G)\langle x\rangle^{-s-1/2-\epsilon}f
\right\|_{L^2}d\tau$$
 $$\le C\gamma\int_{-\infty}^\infty\langle\tau/h\rangle^{-1-\epsilon}
\langle(t-\tau)/h\rangle^{2s}\left\|\langle x\rangle^{-n/2-
\epsilon}e^{i\tau G_0}
\psi_1(h^2G_0)\langle x\rangle^{-s-\epsilon}g\right\|_{L^2}^2d\tau$$
 $$+ C\gamma^{-1}\int_{-\infty}^\infty
\langle\tau/h\rangle^{1+\epsilon}
\left\|\langle x\rangle^{-1-\epsilon}
e^{i\tau G}\psi(h^2G)\langle x\rangle^{-s-1/2-\epsilon}f
\right\|_{L^2}^2d\tau$$
 $$+C\gamma\int_{-\infty}^\infty
\left\|\langle x\rangle^{-1/2-\epsilon}e^{i\tau G_0}
\psi_1(h^2G_0)\langle x\rangle^{-s-\epsilon}g\right\|_{L^2}^2d\tau$$
 $$+C\gamma^{-1}\int_{-\infty}^\infty
\langle\tau/h\rangle^{2s}
\left\|\langle x\rangle^{-(n+1)/2-\epsilon}
e^{i\tau G}\psi(h^2G)\langle x\rangle^{-s-1/2-\epsilon}f
\right\|_{L^2}^2d\tau$$
 $$\le Ch\gamma\|g\|_{L^2}^2+
  Ch\gamma^{-1}\|f\|_{L^2}^2=2Ch\|f\|_{L^2}\|g\|_{L^2},\eqno{(3.18)}$$
if we choose $\gamma=\|f\|_{L^2}/\|g\|_{L^2}$.
By (3.17) and (3.18), we have (with $s=n/2$)
$$\left\|\langle x\rangle^{-s-\epsilon}\Psi(t;h)
\langle x\rangle^{-s-1/2-\epsilon}f\right\|_{L^2}$$ $$\le O(h^2)
\left\|\langle x\rangle^{-(n+2)/2-\epsilon}\Psi(t;h)
\langle x\rangle^{-s-1/2-\epsilon}f\right\|_{L^2}+O(h)
\langle t/h\rangle^{-s}\|f\|_{L^2}.\eqno{(3.19)}$$
Hence, there exists a constant $0<h_0<1$ so that if $0<h\le h_0$,
we can absorbe the first term in the RHS of (3.19), thus obtaining
the estimate (for $0<h\le h_0$)
$$\left\|\langle x\rangle^{-s-\epsilon}\Psi(t;h)
\langle x\rangle^{-s-1/2-\epsilon}\right\|_{L^2\to L^2}\le Ch
\langle t/h\rangle^{-s}.\eqno{(3.20)}$$
Let now $h_0\le h\le 1$. Without loss of generality we may suppose 
$h=1$. By (2.2) and (3.15), the norm in the first term in the RHS
of (3.19) is upper bounded by $C\langle t\rangle^{-s}\|f\|_{L^2}$,
which again implies (3.20). By (2.2) and (3.20), we conclude
$$\left\|\langle x\rangle^{-s-\epsilon}e^{itG}\psi(h^2G)
\langle x\rangle^{-1/2-s-\epsilon}\right\|_{L^2\to L^2}\le C
\langle t/h\rangle^{-s},\eqno{(3.21)}$$
with $s=n/2$, and hence with all $0\le s\le n/2$. To show that this
implies (3.14) with $s=n/2$, we will proceed in the same way as
in Section 3 of \cite{kn:V2}. Let $r=|x|$ denote the radial variable
and set ${\cal D}_r=\langle r\rangle^{-1}rh\partial_r$. It is easy to
see that (3.21) implies
$$\left\|\langle x\rangle^{-s-\epsilon}{\cal D}_re^{itG}\psi(h^2G)
\langle x\rangle^{-1/2-s-\epsilon}\right\|_{L^2\to L^2}\le C
\langle t/h\rangle^{-s},\eqno{(3.22)}$$
for all $0\le s\le n/2$. Furthemore, using Duhamel's formula
together with the identity
$$-2\Delta+[r\partial_r,\Delta]=0,$$
we obtain
$$\psi_1(h^2G)[r\partial_r,e^{itG}]\psi(h^2G)=\int_0^t
\psi_1(h^2G)e^{i(t-\tau)G}[r\partial_r,G]e^{i\tau G}\psi(h^2G)d\tau$$
 $$=\int_0^t
\psi_1(h^2G)e^{i(t-\tau)G}\left(2G-2V+r\partial_rV-Vr\partial_r
\right)e^{i\tau G}\psi(h^2G)d\tau$$
 $$=2te^{itG}G\psi(h^2G)+\int_0^t
\psi_1(h^2G)e^{i(t-\tau)G}\left(-V+\partial_rrV-Vr\partial_r
\right)e^{i\tau G}\psi(h^2G)d\tau,$$
where the functions $\psi$ and $\psi_1$ are as above. Set 
$\widetilde\psi_1(\sigma)=\sigma^{-1}\psi_1(\sigma)$. From the above identity
we get (with $s=n/2$)
$$2(t/h)\langle x\rangle^{-s-\epsilon}e^{itG}\psi(h^2G)
\langle x\rangle^{-s-\epsilon}$$ 
 $$=\left(\langle x\rangle^{-s-\epsilon}\widetilde\psi_1(h^2G)
\langle x\rangle^{s+\epsilon}\right)\langle x\rangle^{-s+1-\epsilon}
{\cal D}_re^{itG}\psi(h^2G)
\langle x\rangle^{-s-\epsilon}$$ 
 $$+\langle x\rangle^{-s-\epsilon}\widetilde\psi_1(h^2G)
e^{itG}\left({\cal D}_r^*+O(h)\right)
\langle x\rangle^{-s+1-\epsilon}
\left(\langle x\rangle^{s+\epsilon}\psi(h^2G)
\langle x\rangle^{-s-\epsilon}\right)$$ 
 $$+h\int_0^\infty\langle x\rangle^{-s-\epsilon}\widetilde
\psi_1(h^2G)e^{i(t-\tau)G}Ve^{i\tau G}\psi(h^2G)
\langle x\rangle^{-s-\epsilon}d\tau$$
 $$+\int_0^\infty\langle x\rangle^{-s-\epsilon}\widetilde
\psi_1(h^2G)e^{i(t-\tau)G}\left({\cal D}_r^*+O(h)\right)
\langle r\rangle Ve^{i\tau G}\psi(h^2G)
\langle x\rangle^{-s-\epsilon}d\tau$$
 $$+\int_0^\infty\langle x\rangle^{-s-\epsilon}\widetilde
\psi_1(h^2G)e^{i(t-\tau)G}V\langle r\rangle{\cal D}_r
e^{i\tau G}\psi(h^2G)\langle x\rangle^{-s-\epsilon}d\tau.
\eqno{(3.23)}$$
By (2.20), (3.21) and (3.22), we have that the $L^2\to L^2$ norm
of each of the first two terms in the RHS of (3.23) is upper
bounded by $O(\langle t/h\rangle^{-s+1})$. Similarly, in view of 
(1.1), we also have that the $L^2\to L^2$ norm
of each integral in the RHS of (3.23) is upper
bounded by $O(\langle t/h\rangle^{-s+1/2})$. Therefore, (3.14)
with $s=n/2$ follows from (3.23) together with the estimates
(3.21) and (3.22).
\eproof

{\it Proof of Proposition 3.4.} We will derive (3.15) from the following 
lemma which can be proved in precisely the same way as Lemma 3.6
of \cite{kn:V2}.

\begin{lemma} Assume (1.1) fulfilled and let $0\le s\le n/2$.  
Let also $m\ge 0$ denote the bigest integer $\le s$ and set
$\mu=s-m$. Then, the
operator-valued function
$${\cal R}_s^\pm(\lambda)=\lambda \langle x\rangle^{-1/2-s-\epsilon}
R^\pm(\lambda)\langle x\rangle^{-1/2-s-\epsilon}:L^2\to L^2$$
is $C^m$ in $\lambda$ for $\lambda>0$, 
$\partial_\lambda^m{\cal R}_s^\pm$ is H\"older of order $\mu$, and
satisfies the estimates
$$\left\|\partial_\lambda^j{\cal R}_s^\pm(\lambda)
\right\|_{L^2\to L^2}\le C,\quad \lambda\ge\lambda_0,\,0\le j\le m,
\eqno{(3.24)}$$
$$\left\|\partial_\lambda^m{\cal R}_s^\pm(\lambda_2)-
\partial_\lambda^m{\cal R}_s^\pm(\lambda_1)
\right\|_{L^2\to L^2}\le C|\lambda_2-\lambda_1|^\mu,\quad 
\lambda_2>\lambda_1\ge\lambda_0,\eqno{(3.25)}$$
for every $\lambda_0>0$, 
with a constant $C>0$ independent of $\lambda$, 
$\lambda_1$ and $\lambda_2$. 
\end{lemma}

We are going to take advantage of the formula (3.8) with 
$h=1$ and $\psi^2$ replaced by $\psi$. Set 
$$T(\lambda)=T^+(\lambda)-T^-(\lambda),\quad 
T^\pm(\lambda)=(2\pi i)^{-1}\lambda^{-1/2}
{\cal R}_s^\pm(\sqrt{\lambda}),\eqno{(3.26)}$$ 
and choose a real-valued
function $\phi\in C_0^\infty([1/3,1/2])$, $\phi\ge 0$, such that 
$\int\phi(\sigma)d\sigma=1$. Then the function
$$T^\pm_\theta(\lambda)=\theta^{-1}\int T^\pm(\lambda+\sigma)\phi(\sigma/
\theta)d\sigma,\quad 0<\theta\le1,$$
is $C^\infty$ in $\lambda$ with values in ${\cal L}(L^2)$ and, 
in view of (3.24) and (3.25), satisfies the estimates  
$$\|\partial_\lambda^jT^\pm_\theta(\lambda)\|_{L^2\to L^2}\le 
C\lambda^{-(j+1)/2},\quad 0\le j\le m,\eqno{(3.27)}$$
 $$\|\partial_\lambda^mT^\pm_\theta(\lambda)-\partial_\lambda^m
T^\pm(\lambda)\|_{L^2\to L^2}\le \theta^{-1}
\int\|\partial_\lambda^mT^\pm(\lambda+\sigma)-\partial_\lambda^m
T^\pm(\lambda)\|_{L^2\to L^2}\phi(\sigma/\theta)d\sigma$$
 $$\le C\lambda^{-(m+1)/2}
\theta^{-1}\int\sigma^{\mu}\phi(\sigma/\theta)d\sigma\le 
 C\lambda^{-(m+1)/2}\theta^{\mu},\eqno{(3.28)}$$ 
$$\|\partial_\lambda^jT^\pm_\theta(\lambda)-\partial_\lambda^j
T^\pm(\lambda)\|_{L^2\to L^2}\le C\lambda^{-(j+1)/2}
\theta,\quad 0\le j\le m-1,
\eqno{(3.29)}$$
 $$\left\|\partial_\lambda^{m+1}T^\pm_\theta(\lambda)
\right\|_{L^2\to L^2}\le 
\theta^{-2}\int\|\partial_\lambda^mT^\pm(\lambda+\sigma)-
\partial_\lambda^mT^\pm(\lambda)\|_{L^2\to L^2}
|\phi'(\sigma/\theta)|d\sigma$$ $$
\le C\lambda^{-(m+2)/2}\theta^{-2}\int\sigma^{\mu}|\phi'(\sigma/\theta)|
d\sigma\le C\lambda^{-(m+2)/2}\theta^{-1+\mu}.\eqno{(3.30)}$$ 
Integrating by parts $m$ times and using 
(3.28) and (3.29), we get
$$\left\|\int_0^\infty
 e^{it\lambda}\psi(\lambda)\left(T^\pm_\theta(\lambda)-
T^\pm(\lambda)\right) d\lambda\right\|_{L^2\to L^2}$$ 
 $$=\left\|t^{-m}\int_0^\infty
 e^{it\lambda}\frac{d^m}{d\lambda^m}\left(
\psi(\lambda)\left(T^\pm_\theta(\lambda)-
T^\pm(\lambda)\right)\right) d\lambda\right\|_{L^2\to L^2}
\le C\theta^{\mu}|t|^{-m}.\eqno{(3.31)}$$
Similarly, integrating by parts $m+1$ times and using (3.27) and 
(3.30), we get
$$\left\|\int_0^\infty e^{it\lambda}\psi(\lambda)
T^\pm_\theta(\lambda) d\lambda\right\|_{L^2\to L^2}$$
 $$=\left\|t^{-m-1}\int_0^\infty
 e^{it\lambda}\frac{d^{m+1}}{d\lambda^{m+1}}\left(\psi(\lambda)
T^\pm_\theta(\lambda)\right)d\lambda\right\|_{L^2\to L^2}
\le C\theta^{-1+\mu}|t|^{-m-1}.\eqno{(3.32)}$$
By (3.26), (3.31) and (3.32), 
$$\left\|\int_0^\infty e^{it\lambda}\psi(\lambda)
T(\lambda) d\lambda\right\|_{L^2\to L^2}\le 
 C\theta^{\mu}|t|^{-m}\left(1+|t|^{-1}\theta^{-1}
\right)\le C|t|^{-m-\mu},\eqno{(3.33)}$$
if we take $\theta=|t|^{-1}$, which clearly implies (3.15).\\

In what follows in this section we will derive (3.16) 
from Lemma 3.5. Let $0\le s\le n/2$ and let $m\ge 0$ be the
bigest integer $\le s$. Remark that the function 
$\partial_\lambda^m{\cal R}_s^\pm$ satisfies (3.25) with
$\mu=s-m+\epsilon/2$. Consequently, the estimates (3.28)
and (3.30) are valid with $\mu=s-m+\epsilon/2$. Let 
$\phi_+\in C^\infty({\bf R})$, $\phi_+(t)=0$ for $t\le 1$,
$\phi_+(t)=1$ for $t\ge 2$. Given any function $f\in L^2$, set
$$u(t;h)=\langle x\rangle^{-1/2-s-\epsilon}\phi_+(t/h)e^{itG}
\varphi(h^2G)\langle x\rangle^{-1/2-s-\epsilon}f.$$
We have
$$\left(\partial_t-iG\right)\langle x\rangle^{1/2+s+\epsilon}u(t;h)$$ $$=
h^{-1}\phi'_+(t/h)e^{itG}
\psi(h^2G)\langle x\rangle^{-1/2-s-\epsilon}f=:h^{-1}
\langle x\rangle^{-1/2-s-\epsilon}v(t;h).\eqno{(3.34)}$$
Clearly,
the support of the function $v(t;h)$ with respect to the variable $t$
is contained in the interval $[h,2h]$, and by (2.20) we have 
$$\|v(t;h)\|_{L^2}\le C\|f\|_{L^2},\quad 1\le t/h\le 2,\eqno{(3.35)}$$
with a constant $C>0$ independent of $t$, $h$ and $f$. Using
Duhamel's formula we deduce from (3.34),  
$$u(t;h)=h^{-1}\int_0^t\langle x\rangle^{-1/2-s-\epsilon}
\psi_1(h^2G)e^{i(t-\tau)G}\langle x\rangle^{-1/2-s-\epsilon}
v(\tau;h)d\tau,\eqno{(3.36)}$$
where the function $\psi_1$ is as above. 
It follows from (3.36) that the Fourier transforms of the functions $u(t;h)$
and $v(t;h)$ satisfy the identity
$$\hat u(\lambda;h)=h^{-1}Q^+(\lambda)\hat v(\lambda;h),
\eqno{(3.37)}$$
where $Q^+(\lambda)$ is the Fourier transform of the operator
$$\langle x\rangle^{-1/2-s-\epsilon}
\psi_1(h^2G)\eta_+(t)e^{itG}
\langle x\rangle^{-1/2-s-\epsilon},$$
$\eta_+$ being the characteristic function of the interval $[0,+\infty)$. 
It is easy to see that
$$Q^+(\lambda)=(2\pi i)^{-1}\langle x\rangle^{-1/2-s-\epsilon}
 \psi_1(h^2G)R^+(\sqrt{\lambda})\langle x\rangle^{-1/2-s-\epsilon}
={\cal B}(h)T^+(\lambda),\eqno{(3.38)}$$
where the operator
$${\cal B}(h)=\langle x\rangle^{-1/2-s-\epsilon}
\psi_1(h^2G)\langle x\rangle^{1/2+s+\epsilon}:L^2\to L^2$$
is bounded uniformly in $h$ in view of (2.20). 
Fix a constant $0<\gamma<1$ such that supp$\,\psi_1\subset(\gamma,
\gamma^{-1})$. Then , for $\lambda h^2\in{\bf R}\setminus
(\gamma,\gamma^{-1})$, we have
$$\left\|\frac{d^kQ^+(\lambda)}{d\lambda^k}\right\|_{L^2\to L^2}
\le C\left\|\psi_1(h^2G)(G-\lambda)^{-k-1}\right\|_{L^2\to L^2}$$
 $$\le C\sup_{\sigma\in{\bf R}}\left|\psi_1(h^2\sigma)(\sigma-\lambda)^{-k-1}
\right|\le C_kh^{2k+2},\eqno{(3.39)}$$
for every integer $k\ge 0$ with a constant $C_k>0$ independent
of $\lambda$ and $h$. 
Set $Q^+_\theta(\lambda)={\cal B}(h)T^+_\theta(\lambda)$ 
and define the function $u_\theta(t;h)$ via the formula
$$\hat u_\theta(\lambda;h)=Q^+_\theta(\lambda)\hat v(\lambda;h).$$
Using (3.27)-(3.30) when $\lambda h^2\in(\gamma,\gamma^{-1})$
and (3.39) when $\lambda h^2\in{\bf R}\setminus
(\gamma,\gamma^{-1})$, we obtain
$$\|\partial_\lambda^jQ^+_\theta(\lambda)\|_{L^2\to L^2}\le 
Ch^{j+1},\quad 0\le j\le m,\eqno{(3.40)}$$
 $$\|\partial_\lambda^mQ^+_\theta(\lambda)-\partial_\lambda^m
Q^+(\lambda)\|_{L^2\to L^2}\le \theta^{-1}
\int\|\partial_\lambda^mQ^+(\lambda+\sigma)-\partial_\lambda^m
Q^+(\lambda)\|_{L^2\to L^2}\phi(\sigma/\theta)d\sigma$$
 $$\le Ch^{m+1}
\theta^{-1}\int\sigma^{\mu}\phi(\sigma/\theta)d\sigma\le 
 Ch^{m+1}\theta^{\mu},\eqno{(3.41)}$$ 
$$\|\partial_\lambda^jQ^+_\theta(\lambda)-\partial_\lambda^j
Q^+(\lambda)\|_{L^2\to L^2}\le Ch^{j+1}
\theta,\quad 0\le j\le m-1,\eqno{(3.42)}$$
 $$\left\|\partial_\lambda^{m+1}Q^+_\theta(\lambda)
\right\|_{L^2\to L^2}\le 
\theta^{-2}\int\|\partial_\lambda^mQ^+(\lambda+\sigma)-
\partial_\lambda^mQ^+(\lambda)\|_{L^2\to L^2}
|\phi'(\sigma/\theta)|d\sigma$$ $$
\le Ch^{m+2}\theta^{-2}\int\sigma^{\mu}|\phi'(\sigma/\theta)|
d\sigma\le Ch^{m+2}\theta^{-1+\mu}.\eqno{(3.43)}$$ 
Using (3.40)-(3.43) (with $\mu=s-m+\epsilon/2$) together with the 
Plancherel identity and (3.35), we obtain
$$\int_{-\infty}^\infty|t|^{2m}\|u_\theta(t;h)-u(t;h)\|_{L^2}^2dt
=\int_{-\infty}^\infty\|\partial_\lambda^m(\hat u_\theta(\lambda;h)
-\hat u(\lambda;h))\|_{L^2}^2d\lambda$$ 
 $$\le Ch^{-2}\sum_{k=0}^m\int_{-\infty}^\infty\left\|\left(
\partial_\lambda^kQ^+_\theta(\lambda)-\partial_\lambda^k
Q^+(\lambda)\right)
\partial_\lambda^{m-k}\hat v(\lambda;h)\right\|_{L^2}^2d\lambda$$ 
 $$\le Ch^{2m}\theta^{2\mu}\sum_{k=0}^m\int_{-\infty}^\infty\left\|
\partial_\lambda^{m-k}\hat v(\lambda;h)\right\|_{L^2}^2d\lambda
=Ch^{2m}\theta^{2\mu}\sum_{k=0}^m\int_{-\infty}^\infty
|t|^{2m-2k}\|v(t;h)\|_{L^2}^2dt$$
 $$\le Ch^{2m}\theta^{2\mu}\int_h^{2h}\|v(t;h)\|_{L^2}^2dt
\le  Ch^{2m+1}\theta^{2\mu}\|f\|_{L^2}^2,\eqno{(3.44)}$$
with a constant $C>0$ independent of $h$, $\theta$ and $f$.
By (3.44) we get, $\forall A\ge 1$,
$$\int_{Ah}^{2Ah}\|u_\theta(t;h)-u(t;h)\|_{L^2}^2dt
\le  ChA^{-2m}\theta^{2\mu}\|f\|_{L^2}^2.\eqno{(3.45)}$$
In the same way, we obtain
$$\int_{-\infty}^\infty|t|^{2m+2}\|u_\theta(t;h)\|_{L^2}^2dt
=\int_{-\infty}^\infty\|\partial_\lambda^{m+1}\hat u_\theta(\lambda;h)
\|_{L^2}^2d\lambda$$ 
 $$\le Ch^{-2}\sum_{k=0}^{m+1}\int_{-\infty}^\infty\left\|
\partial_\lambda^kQ^+_\theta(\lambda)
\partial_\lambda^{m+1-k}\hat v(\lambda;h)\right\|_{L^2}^2d\lambda$$ 
 $$\le Ch^{2m+2}\theta^{-2+2\mu}\sum_{k=0}^{m+1}\int_{-\infty}^\infty\left\|
\partial_\lambda^{m+1-k}\hat v(\lambda;h)\right\|_{L^2}^2d\lambda$$ $$
=Ch^{2m+2}\theta^{-2+2\mu}\sum_{k=0}^{m+1}\int_{-\infty}^\infty
|t|^{2m+2-2k}\|v(t;h)\|_{L^2}^2dt$$
 $$\le Ch^{2m+2}\theta^{-2+2\mu}\int_h^{2h}\|v(t;h)\|_{L^2}^2dt
\le  Ch^{2m+3}\theta^{-2+2\mu}\|f\|_{L^2}^2,\eqno{(3.46)}$$
with a constant $C>0$ independent of $h$, $\theta$ and $f$.
By (3.46) we get, $\forall A\ge 1$,
$$\int_{Ah}^{2Ah}\|u_\theta(t;h)\|_{L^2}^2dt
\le  ChA^{-2m-2}\theta^{-2+2\mu}\|f\|_{L^2}^2.\eqno{(3.47)}$$
Combining (3.45) and (3.47) leads to
$$\int_{Ah}^{2Ah}(t/h)^{2s}\|u(t;h)\|_{L^2}^2dt
\le  ChA^{2\mu-\epsilon}\theta^{2\mu}
\left(1+A^{-2}\theta^{-2}\right)\|f\|_{L^2}^2\le
  ChA^{-\epsilon}\|f\|_{L^2}^2,\eqno{(3.48)}$$
if we choose $\theta=A^{-1}$, where $s=m+\mu-\epsilon/2$.
By (3.48), for every integer $k\ge 0$ we have
$$\int_{2^kh}^{2^{k+1}h}(t/h)^{2s}\|u(t;h)\|_{L^2}^2dt
\le Ch2^{-\epsilon k}\|f\|_{L^2}^2.\eqno{(3.49)}$$
Summing up (3.49) leads to
$$\int_h^{\infty}(t/h)^{2s}\|u(t;h)\|_{L^2}^2dt
\le Ch\|f\|_{L^2}^2,\eqno{(3.50)}$$
which clearly implies (3.16).
\eproof

\section{Proof of Theorem 1.1}
We will first prove the following

\begin{prop} For every $0<\epsilon\ll 1$, $1/2-\epsilon/2\le s\le
(n-1)/2$, $0<h\le 1$, $t\neq 0$, 
 we have 
$$\left\|\Psi(t;h)\langle x\rangle^{-s-1/2-\epsilon}
\right\|_{L^2\to L^\infty}\le Ch^{-(n-3)/2+s-\epsilon}|t|^{-s-1/2},
\eqno{(4.1)}$$
with a constant $C>0$ independent of $t$ and $h$.
\end{prop}

{\it Proof.} By (2.1) and (2.24), we have
$$\left\|\Psi_1(t;h)\langle x\rangle^{-s-1/2-\epsilon}f
\right\|_{L^\infty}$$ $$\le O(h^2)
\left\|\Psi(t;h)\langle x\rangle^{-s-1/2-\epsilon}f
\right\|_{L^\infty}+
Ch^{-(n-3)/2+s}|t|^{-s-1/2}\|f\|_{L^2},\quad 
\forall f\in L^2.\eqno{(4.2)}$$
Using (2.1) and (3.14), we obtain
$$\left\|\Psi_2(t;h)\langle x\rangle^{-s-1/2-\epsilon}
\right\|_{L^2\to L^\infty}$$ $$\le
C\int_0^{t/2}\left\|e^{i(t-\tau)G_0}\psi_1(h^2G_0)
\langle x\rangle^{-s-1/2-\epsilon}\right\|_{L^2\to L^\infty}
\left\|\langle x\rangle^{-1-\epsilon}e^{i\tau G}\psi(h^2G)
\langle x\rangle^{-1-\epsilon/2}\right\|_{L^2\to L^2}d\tau$$
 $$+C\int_0^{t/2}\left\|e^{i\tau G_0}\psi_1(h^2G_0)
\langle x\rangle^{-1-\epsilon}\right\|_{L^2\to L^\infty}
\left\|\langle x\rangle^{-1/2-s-\epsilon}e^{i(t-\tau)G}\psi(h^2G)
\langle x\rangle^{-1/2-s-\epsilon}\right\|_{L^2\to L^2}d\tau $$
 $$\le Ch^{s-(n-1)/2}|t|^{-s-1/2}\int_0^{t/2}\langle\tau/h
\rangle^{-1-\epsilon/2}d\tau
+Ch^{-(n-2)/2-\epsilon}\langle t/h\rangle^{-s-1/2}
\int_0^{t/2}|\tau|^{-1+\epsilon/2}\langle\tau
\rangle^{-\epsilon}d\tau$$ $$\le 
Ch^{s-(n-3)/2-\epsilon}|t|^{-s-1/2}.\eqno{(4.3)}$$
By (4.2) and (4.3),
$$\left\|\Psi(t;h)\langle x\rangle^{-s-1/2-\epsilon}f
\right\|_{L^\infty}$$ $$\le O(h^2)
\left\|\Psi(t;h)\langle x\rangle^{-s-1/2-\epsilon}f
\right\|_{L^\infty}+
Ch^{-(n-3)/2+s-\epsilon}|t|^{-s-1/2}\|f\|_{L^2},\quad 
\forall f\in L^2.\eqno{(4.4)}$$
Hence, there exists a constant $0<h_0<1$ so that for 
$0<h\le h_0$ we can absorbe the first term in the RHS
of (4.4), thus obtaining (4.1) in this case. Let now
$h_0\le h\le 1$. Without loss of generality we may suppose
$h=1$. Then the only term we need to estimate is
$$\left\|\left(\psi_1(G)-\psi_1(G_0)\right)e^{itG}\psi(G)
\langle x\rangle^{-s-1/2-\epsilon}f\right\|_{L^\infty}.$$
In view of (2.27), this is reduced to estimating
$$\left\|\langle x\rangle^{-(n+2)/2-\epsilon}e^{itG}\psi(G)
\langle x\rangle^{-s-1/2-\epsilon}f\right\|_{L^2},$$
which, in view of Theorem 3.3, is upper bounded by
$O(|t|^{-s-1/2})\|f\|_{L^2}$.
\eproof

Write $\Psi_2=\Psi_3+\Psi_4$, where
$$\Psi_3(t;h)=i\int_0^t e^{i(t-\tau)G_0}\psi_1(h^2G_0)Ve^{i
\tau G_0}\psi(h^2G_0)d\tau,$$
$$\Psi_4(t;h)=i\int_0^t e^{i(t-\tau)G_0}\psi_1(h^2G_0)V\Psi(\tau;h)d\tau.$$

\begin{prop} For every $0<\epsilon\ll 1$, $0<h\le 1$, $t\neq 0$, 
 we have  
$$\left\|\Psi_1(t;h)\right\|_{L^1\to L^\infty}
\le Ch^{2-n/2}|t|^{-n/2},\eqno{(4.5)}$$
$$\left\|\Psi_4(t;h)\right\|_{L^1\to L^\infty}
\le Ch^{2-n/2-\epsilon}|t|^{-n/2},\eqno{(4.6)}$$ 
with a constant $C>0$ independent of $t$ and $h$.
\end{prop}

{\it Proof.} By (2.24), (2.27) and (4.1), we have
$$\left\|\Psi_1(t;h)\right\|_{L^1\to L^\infty}
\le Ch^{2-n/2}\left\|\langle x\rangle^{-(n+2)/2-\epsilon}
\Psi(t;h)\right\|_{L^1\to L^2}$$ $$+Ch^2|t|^{-n/2}
\le Ch^{3-n/2-\epsilon}|t|^{-n/2},$$
which implies (4.5). By (2.1) and (4.1), we have
$$\left\|\Psi_4(t;h)\right\|_{L^1\to L^\infty}$$ $$
\le C\int_0^{t/2}\left\|e^{i(t-\tau)G_0}\psi_1(h^2G_0)
\langle x\rangle^{-n/2-\epsilon}\right\|_{L^2\to L^\infty}
\left\|\langle x\rangle^{-1-\epsilon}e^{i\tau G}\psi(h^2G)
\right\|_{L^1\to L^2}d\tau$$
 $$+ C\int_0^{t/2}\left\|e^{i\tau G_0}\psi_1(h^2G_0)
\langle x\rangle^{-1-\epsilon}\right\|_{L^2\to L^\infty}
\left\|\langle x\rangle^{-n/2-\epsilon}e^{i(t-\tau)G}\psi(h^2G)
\right\|_{L^1\to L^2}d\tau$$
  $$\le 2Ch^{-(n-4)/2-\epsilon}|t|^{-n/2}\int_0^{t/2}
|\tau|^{-1+\epsilon/2}\langle\tau
\rangle^{-\epsilon}d\tau\le C_\epsilon 
h^{-(n-4)/2-\epsilon}|t|^{-n/2},$$
 $\forall 0<\epsilon\ll 1$, with a constant $C_\epsilon>0$
independent of $t$ and $h$.
\eproof

We will now derive Theorem 1.1 from the estimates (4.1), (4.5), (4.6)
and the following

\begin{prop}  For every $0<h\le 1$, $t\neq 0$, 
 we have  
$$\left\|\Psi_3(t;h)\right\|_{L^1\to L^\infty}
\le Ch^{-(n-3)/2}|t|^{-n/2},\eqno{(4.7)}$$ 
with a constant $C>0$ independent of $t$ and $h$.
Moreover, the operator $\Psi_3$ is of the form
$$\Psi_3(t;h)=E(t;h)+F(t;h),
\eqno{(4.8)}$$
where the operator $E$ has a kernel of the form
$$\int_{{\bf R}^n}\int w(\lambda,t,|x-\xi|,|y-\xi|) 
\lambda^{(n-3)/4}\psi(h^2\lambda)
V(\xi)d\lambda d\xi,\eqno{(4.9)}$$
with a function $w$ independent of $h$ and satisfying the bound, 
$\forall t\neq 0,\,\sigma_1,\sigma_2>0$,
$$\left|\int w(\lambda,t,\sigma_1,\sigma_2)\chi_a(\lambda)d\lambda
\right|\le C|t|^{-n/2}\left(\sigma_1^{-(n-1)/2}+\sigma_1^{-n+1}+
\sigma_2^{-(n-1)/2}+\sigma_2^{-n+1}
\right).\eqno{(4.10)}$$ 
The operator $F$ satisfies
$$\left\|F(t;h)\right\|_{L^1\to L^\infty}
\le Ch^{2-n/2}|t|^{-n/2},\quad\forall t\neq 0,\,0<h\le 1,\eqno{(4.11)}$$ 
with a constant $C>0$ independent of $t$ and $h$.
\end{prop}

Writting the function $\chi_a$ as
$$\chi_a(\sigma)=\int_0^1\psi(\sigma\theta)\frac{d\theta}{\theta},$$
where $\psi(\sigma)=\sigma\chi'_a(\sigma)\in C_0^\infty((0,+\infty))$,
we obtain by (4.1) with $s=(n-1)/2$, $0<\epsilon\ll 1$, 
$$\left\|e^{itG}\chi_a(G)\langle x\rangle^{-n/2-\epsilon}-
e^{itG_0}\chi_a(G_0)\langle x\rangle^{-n/2-\epsilon}\right\|_{L^2\to 
L^\infty}$$
 $$\le \int_0^1\left\|\Psi(t;\sqrt{\theta})\langle x
\rangle^{-n/2-\epsilon}\right\|_{L^2\to 
L^\infty}\frac{d\theta}{\theta}\le C|t|^{-n/2}\int_0^1
\theta^{-1/2-\epsilon/2}d\theta\le C|t|^{-n/2},\eqno{(4.12)}$$
which implies (1.3). 

Take $\psi(\sigma)=\sigma^{1-(n-3)/4}\chi'_a(\sigma)$ and 
denote by ${\cal E}(t)$ the operator with kernel defined by
replacing in (4.9) the function $\lambda^{(n-3)/4}\psi(h^2\lambda)$ by 
$\chi_a(\lambda)$.
By (1.1) and (4.10), we have
$$\left\|{\cal E}(t)\right\|_{L^1\to L^\infty}
\le C|t|^{-n/2},\quad\forall t\neq 0.\eqno{(4.13)}$$ 
By (4.5), (4.6) and (4.11) 
$$\left\|G^{-(n-3)/4}e^{itG}\chi_a(G)-
G_0^{-(n-3)/4}e^{itG_0}\chi_a(G_0)-{\cal E}(t)
\right\|_{L^1\to L^\infty}$$
 $$\le \int_0^1\left\|\Psi(t;\sqrt{\theta})-E(t;\sqrt{\theta})
\right\|_{L^1\to 
L^\infty}\theta^{-1+(n-3)/4}d\theta$$ $$\le C|t|^{-n/2}\int_0^1
\theta^{-3/4-\epsilon/2}d\theta\le C|t|^{-n/2},\eqno{(4.14)}$$
which together with (4.13) and the fact that the operators 
$G_0^{-(n-3)/4}\chi_a(G_0)$ and $\chi_a(G_0)$ are bounded 
on $L^p$, $1\le p\le\infty$, imply (1.2). 

To prove (1.4), observe
that by (4.5), (4.6) and (4.7) we have
$$\left\|\Psi(t;h)\right\|_{L^1\to L^\infty}
\le Ch^{-(n-3)/2}|t|^{-n/2},\quad \forall t\neq 0,\,0<h\le 1.\eqno{(4.15)}$$ 
By interpolation between (3.1) and (4.15), we get
$$\left\|\Psi(t;h)\right\|_{L^{p'}\to L^p}
\le Ch^{1-\alpha(n-1)/2}|t|^{-\alpha n/2},
\quad \forall t\neq 0,\,0<h\le 1,\eqno{(4.16)}$$ 
for every $2\le p\le +\infty$, where $1/p+1/p'=1$, $\alpha=1-2/p$.
As above, taking $\psi(\sigma)=\sigma^{1-\alpha q/2}\chi'_a(\sigma)$ and using
(4.16), we get
$$\left\|G^{-\alpha q/2}e^{itG}\chi_a(G)-G_0^{-\alpha q/2}
e^{itG_0}\chi_a(G_0)\right\|_{L^{p'}\to L^p}
 \le \int_0^1\left\|\Psi(t;\sqrt{\theta})\right\|_{L^{p'}\to 
L^p}\theta^{-1+\alpha q/2}d\theta$$ $$\le C|t|^{-\alpha n/2}\int_0^1
\theta^{-1/2-\alpha(n-1)/4+\alpha q/2}
d\theta\le C|t|^{-\alpha n/2}\eqno{(4.17)}$$
if $1/2+\alpha(n-1)/4-\alpha q/2<1$, that is, for $2\le p<2(n-1-2q)/(n-3-2q)$. 
This clearly proves (1.4). 
\eproof

\section{Proof of Proposition 4.3} 

The kernel of the operator $\Psi_3$ is of the form
$$\int_{{\bf R}^n}U_h(|x-\xi|,|y-\xi|;t)V(\xi)d\xi,\eqno{(5.1)}$$
where
$$U_h(\sigma_1,\sigma_2;t)=i\int_0^t\widetilde K_h(\sigma_1,t-\tau)
K_h(\sigma_2,\tau)d\tau =h^{-2n+2}U_1(\sigma_1h^{-1},\sigma_2h^{-1};
th^{-2}),\eqno{(5.2)}$$
where $K_h$ and $\widetilde K_h$ are defined by (2.4) (in the case
of $\widetilde K_h$ the function $\psi$ is replaced by $\psi_1$).
It is easy to see that to prove the proposition it suffices to show
that the function $U_1$ can be decomposed as $U_1=W_1+L_1$ with
a function $W_1$ of the form
$$W_1(\sigma_1,\sigma_2;t)=\int w(\lambda,t,\sigma_1,\sigma_2) 
\lambda^{(n-3)/4}\psi(\lambda)d\lambda,\eqno{(5.3)}$$
with a function $w$ independent of $\psi$, satisfying (4.10) and
$$w(h^2\lambda,th^{-2},\sigma_1h^{-1},\sigma_2h^{-1})=
h^{(3n-5)/2}w(\lambda,t,\sigma_1,\sigma_2),\eqno{(5.4)}$$ 
while the function $L_1$ satisfies
$$|L_1(\sigma_1,\sigma_2;t)|\le C|t|^{-n/2}
\left(\frac{
\langle\sigma_1\rangle^{(n-2)/2}}{\sigma_1^{n-1}}
+\frac{\langle\sigma_2\rangle^{(n-2)/2}}{\sigma_2^{n-1}}
\right),\quad\forall t\neq 0,\,\sigma_1,
\sigma_2>0.\eqno{(5.5)}$$
To do so, observe that the function $U_1$
is of the form $U_1=U_1^{(1)}+U_1^{(2)}$, where 
$$U_1^{(1)}(\sigma_1,\sigma_2;t)=\frac{(\sigma_1\sigma_2)^{-2\nu}}
{(2\pi)^n}\int\int e^{it\lambda_1^2}\psi_1(\lambda_1^2)\psi(\lambda_2^2)
{\cal J}_\nu(\sigma_1\lambda_1){\cal J}_\nu(\sigma_2\lambda_2)
\frac{\lambda_1\lambda_2}{\lambda_1^2-\lambda_2^2}d\lambda_1d\lambda_2,$$
 $$U_1^{(2)}(\sigma_1,\sigma_2;t)=\frac{(\sigma_1\sigma_2)^{-2\nu}}
{(2\pi)^n}\int\int e^{it\lambda_2^2}\psi_1(\lambda_1^2)\psi(\lambda_2^2)
{\cal J}_\nu(\sigma_1\lambda_1){\cal J}_\nu(\sigma_2\lambda_2)
\frac{\lambda_1\lambda_2}{\lambda_2^2-\lambda_1^2}d\lambda_1d\lambda_2.$$
Recall that ${\cal J}_\nu(z)=e^{iz}b_\nu^+(z)+e^{-iz}b_\nu^-(z)$
with functions $b_\nu^\pm$ satisfying (2.11). Set 
${\cal N}_\nu(z)=e^{iz}b_\nu^+(z)-e^{-iz}b_\nu^-(z)$, and
$$a^\pm(\lambda_1,\lambda_2;\sigma_2)=(\lambda_1-\lambda_2)^{-1}
\left(\frac{\lambda_2}{\lambda_1+\lambda_2}\psi(\lambda_2^2)b_\nu^\pm(
\sigma_2\lambda_2)-\frac{1}{2}\psi(\lambda_1^2)b_\nu^\pm(
\sigma_2\lambda_1)\right).$$
We have 
$$\int_{-\infty}^\infty \psi(\lambda_2^2){\cal J}_\nu(\sigma_2\lambda_2)
\frac{\lambda_2d\lambda_2}{\lambda_1^2-\lambda_2^2}
=\sum_{\pm}\int_{-\infty}^\infty e^{\pm i\sigma_2\lambda_2}
a^\pm(\lambda_1,\lambda_2;\sigma_2)d\lambda_2$$
 $$+\sum_{\pm}\frac{1}{2}\psi(\lambda_1^2)b_\nu^\pm(
\sigma_2\lambda_1)\int_{-\infty}^\infty e^{\pm i\sigma_2\lambda_2}
\frac{d\lambda_2}{\lambda_1-\lambda_2}$$
 $$=\sum_{\pm}\int_{-\infty}^\infty e^{\pm i\sigma_2\lambda_2}
a^\pm(\lambda_1,\lambda_2;\sigma_2)d\lambda_2
 +Const\sum_{\pm}\pm\psi(\lambda_1^2)b_\nu^\pm(
\sigma_2\lambda_1)e^{\pm i\sigma_2\lambda_1}$$
 $$=:\sum_{\pm}A^\pm(\lambda_1,\sigma_2)+Const\,\psi(\lambda_1^2)
{\cal N}_\nu(\sigma_2\lambda_1),\eqno{(5.6)}$$
where we have used that, for any $\sigma>0$, 
$$\int_{-\infty}^\infty e^{\pm i\sigma\mu}\frac{d\mu}{\mu}=
\pm 2i\int_0^\infty\frac{\sin(\sigma\mu)}{\mu}d\mu=
\pm 2i\int_0^\infty\frac{\sin\mu}{\mu}d\mu=\pm 2i\,Const.$$
In view of (5.6) we can write the function $U^{(1)}_1$ in the form
$$U_1^{(1)}(\sigma_1,\sigma_2;t)=Const(\sigma_1\sigma_2)^{-2\nu}
\int e^{it\lambda_1^2}\psi(\lambda_1^2)\
{\cal J}_\nu(\sigma_1\lambda_1){\cal N}_\nu(\sigma_2\lambda_1)
\lambda_1d\lambda_1$$ 
 $$+\frac{(\sigma_1\sigma_2)^{-2\nu}}
{(2\pi)^n}\int e^{it\lambda_1^2}\psi_1(\lambda_1^2)
{\cal J}_\nu(\sigma_1\lambda_1)A(\lambda_1,\sigma_2)
\lambda_1d\lambda_1$$ $$=:W_1^{(1)}
(\sigma_1,\sigma_2;t)+L_1^{(1)}(\sigma_1,\sigma_2;t),
\eqno{(5.7)}$$
where $A=A^++A^-$. 
We will now show that the function $L_1^{(1)}$
satisfies (5.5). Observe first that the functions $a^\pm$
satisfy (for $\lambda_1^2\in\mbox{supp}\,\psi_1$)
$$\left|\partial_{\lambda_1}^{\alpha_1}\partial_{\lambda_2}^{\alpha_2}
a^\pm(\lambda_1,\lambda_2;\sigma)\right|\le C\langle\sigma\rangle^{(n-3)/2}
\langle\lambda_2\rangle^{-1-\alpha_2},\quad\forall\lambda_2\in{\bf R},\,
\sigma>0,\eqno{(5.8)}$$
for all multi-indices $(\alpha_1,\alpha_2)$. Indeed, it is easy to
see that for $\sigma\ge 1$ the bound (5.8) follows from (2.11), while
for $0<\sigma\le 1$ one needs to use the fact that near $z=0$
the functions $b^\pm_\nu$ are of the form
$$b^\pm_\nu(z)=b^\pm_{\nu,1}(z)+z^{n-2}\log zb^\pm_{\nu,2}(z),\eqno{(5.9)}$$
where the functions $b^\pm_{\nu,j}$ are analytic at $z=0$, 
$b^\pm_{\nu,2}\equiv 0$ if $n$ is odd. Therefore, we have 
(for $\lambda^2\in\mbox{supp}\,\psi_1$)
$$\left|\frac{d^k}{d\lambda^k}b^\pm_\nu(\sigma\lambda)\right|\le
C_k,\quad 0<\sigma\le 1,\eqno{(5.10)}$$
for every integer $k$ with a constant $C_k>0$ independent of 
$\sigma$. Clearly, (5.8) for $0<\sigma\le 1$ follows from (5.10). 
Furthermore, an integration by parts together with (5.8) lead to the 
following bounds for the functions $A^\pm$ 
(for $\lambda^2\in\mbox{supp}\,\psi_1$)
$$\left|\partial_\lambda^\alpha A^\pm(\lambda,\sigma)\right|\le C_{\alpha,k}
\langle\sigma\rangle^{(n-3)/2}\sigma^{-k},\quad\forall\sigma>0,
\eqno{(5.11)}$$
for every integers $\alpha\ge 0$, $k\ge 1$. Hence,
$$\left|\partial_\lambda^\alpha A^\pm(\lambda,\sigma)\right|\le C_{\alpha}
\sigma^{-1},\quad\forall\sigma>0.\eqno{(5.12)}$$
On the other hand, by (2.8) and (2.11), we have
$$\left|\partial_z^k{\cal J}_\nu(z)\right|\le Cz^{n-2-k}\langle z
\rangle^{k-(n-1)/2},\quad\forall z>0,\eqno{(5.13)}$$
for every integer $0\le k\le n-2$, while for $k\ge n-1$ we have
$$\left|\partial_z^k{\cal J}_\nu(z)\right|\le C_k\langle z
\rangle^{(n-3)/2},\quad\forall z>0.\eqno{(5.14)}$$
Thus, by (5.12)-(5.14) we obtain (for $\lambda^2\in\mbox{supp}\,\psi_1$)
$$\left|\frac{d^m}{d\lambda^m}\left({\cal J}_\nu(\sigma_1\lambda)A(
\lambda,\sigma_2)\right)\right|\le C\sum_{k=0}^m
\left|\frac{d^k}{d\lambda^k}{\cal J}_\nu(\sigma_1\lambda)
\frac{d^{m-k}}{d\lambda^{m-k}}A(\lambda,\sigma_2)\right|$$
 $$\le C\sigma_2^{-1}\sigma_1^{n-2}\langle\sigma_1\rangle^{m-(n-1)/2},
\eqno{(5.15)}$$
for every integer $m\ge 0$. In the same way as in the proof of (2.5),
using (2.7) together with (5.15), we deduce
$$\left|L_1^{(1)}(\sigma_1,\sigma_2;t)\right|\le C_m|t|^{-m-1/2}
\sigma_2^{-n+1}\langle\sigma_1\rangle^{m-(n-1)/2},
\eqno{(5.16)}$$
for every integer $m\ge 0$, and hence for all real $m\ge 0$. Taking
$m=(n-1)/2$ in (5.16) we obtain the desired bound. Let 
$\phi\in C_0^\infty([-1,1])$, $\phi=1$ on $[-1/2,1/2]$. We further
decompose the function $W_1^{(1)}$ as follows
$$W_1^{(1)}(\sigma_1,\sigma_2;t)=Const(\sigma_1\sigma_2)^{-n+2}
\int e^{it\lambda^2}\psi(\lambda^2)
{\cal J}_\nu(\sigma_1\lambda){\cal N}_\nu(\sigma_2\lambda)
\phi(\sigma_2\lambda)
\lambda d\lambda$$ $$+Const(\sigma_1\sigma_2)^{-n+2}
\int e^{it\lambda^2}\psi(\lambda^2)
{\cal J}_\nu(\sigma_1\lambda)\phi(\sigma_1\lambda)
{\cal N}_\nu(\sigma_2\lambda)(1-\phi)
(\sigma_2\lambda)\lambda d\lambda$$ 
 $$+Const(\sigma_1\sigma_2)^{-n+2}
\int e^{it\lambda^2}\psi(\lambda^2)
{\cal J}_\nu(\sigma_1\lambda)(1-\phi)(\sigma_1\lambda)
{\cal N}_\nu(\sigma_2\lambda)(1-\phi)
(\sigma_2\lambda)\lambda d\lambda$$ 
$$=:M_1^{(1)}(\sigma_1,\sigma_2;t)+N_1^{(1)}(\sigma_1,\sigma_2;t)
+\widetilde W_1^{(1)}(\sigma_1,\sigma_2;t).\eqno{(5.17)}$$
By (5.10) we have (for $\lambda^2\in\mbox{supp}\,\psi$)
$$\left|\frac{d^k}{d\lambda^k}({\cal N}_\nu\phi)(\sigma\lambda)\right|\le
C_k,\quad \forall\sigma>0,\eqno{(5.18)}$$
for every integer $k$ with a constant $C_k>0$ independent of 
$\sigma$. By (5.13), (5.14) and (5.18) we get 
(for $\lambda^2\in\mbox{supp}\,\psi$)
$$\left|\frac{d^m}{d\lambda^m}\left({\cal J}_\nu(\sigma_1\lambda)
({\cal N}_\nu\phi)(\sigma_2\lambda)\right)\right|\le
 C\sigma_1^{n-2}\langle\sigma_1\rangle^{m-(n-1)/2},\eqno{(5.19)}$$
for every integer $m\ge 0$. In the same way as above we deduce from
(5.19)
$$\left|M_1^{(1)}(\sigma_1,\sigma_2;t)\right|\le C_m|t|^{-m-1/2}
\sigma_2^{-n+2}\langle\sigma_1\rangle^{m-(n-1)/2},
\eqno{(5.20)}$$
for all real $m\ge 0$. Taking $m=(n-1)/2$ we conclude that 
$M_1^{(1)}$ satisfies (5.5). Furthermore, in view of (5.13) and (5.14), 
we have (for $\lambda^2\in\mbox{supp}\,\psi$)
$$\left|\frac{d^k}{d\lambda^k}({\cal J}_\nu\phi)(\sigma\lambda)\right|\le
C_k\sigma^{n-2},\quad \forall\sigma>0,\eqno{(5.21)}$$
for every integer $k$ with a constant $C_k>0$ independent of 
$\sigma$, while (2.11) leads to the bound
$$\left|\frac{d^k}{d\lambda^k}((1-\phi){\cal N}_\nu)(\sigma\lambda)\right|\le
C_k\langle\sigma\rangle^{k-(n-1)/2},\quad \forall\sigma>0.\eqno{(5.22)}$$
By (5.21) and (5.22), (for $\lambda^2\in\mbox{supp}\,\psi$)
$$\left|\frac{d^m}{d\lambda^m}\left((\phi{\cal J}_\nu)(\sigma_1\lambda)
((1-\phi){\cal N}_\nu)(\sigma_2\lambda)\right)\right|\le
 C\sigma_1^{n-2}\langle\sigma_2\rangle^{m-(n-1)/2},\eqno{(5.23)}$$
for every integer $m\ge 0$. By (5.23) we get
$$\left|N_1^{(1)}(\sigma_1,\sigma_2;t)\right|\le C_m|t|^{-m-1/2}
\sigma_2^{-n+2}\langle\sigma_2\rangle^{m-(n-1)/2},
\eqno{(5.24)}$$
for all real $m\ge 0$. Taking $m=(n-1)/2$ we conclude that 
$N_1^{(1)}$ satisfies (5.5), too. 

We will now decompose $\widetilde W_1^{(1)}$ using that the functions
$b_\nu^\pm$ admit the expansion
$$b_\nu^\pm(z)=\pm c_\nu z^{(n-3)/2}+O\left(z^{(n-5)/2}\right),\quad
z\to +\infty,$$
where $c_\nu$ is some constant. More precisely, we have
$$\left|\partial_z^k\left(b_\nu^\pm(z)\mp c_\nu z^{(n-3)/2}\right)\right|
\le C_kz^{(n-5)/2-k},\quad z\ge z_0,\eqno{(5.25)}$$
for every integer $k\ge 0$ and every $z_0>0$, with a constant
$C_k>0$ independent of $z$ but depending on $k$ and $z_0$. Thus we can write
$$(\sigma_1\sigma_2)^{-n+2}\lambda((1-\phi)
{\cal J}_\nu)(\sigma_1\lambda)((1-\phi)
{\cal N}_\nu)(\sigma_2\lambda)$$
 $$=c_\nu^2(\sigma_1\sigma_2)^{-(n-1)/2}\lambda^{n-2}
\left(e^{i\sigma_1\lambda}-
e^{-i\sigma_1\lambda}\right)\left(e^{i\sigma_2\lambda}+
e^{-i\sigma_2\lambda}\right)$$
 $$+c_\nu^2(\sigma_1\sigma_2)^{-(n-1)/2}\lambda^{n-2}(-\phi
(\sigma_1\lambda)-\phi(\sigma_2\lambda)+\phi
(\sigma_1\lambda)\phi(\sigma_2\lambda))\left(e^{i\sigma_1\lambda}-
e^{-i\sigma_1\lambda}\right)\left(e^{i\sigma_2\lambda}+
e^{-i\sigma_2\lambda}\right)$$
 $$+(\sigma_1\sigma_2)^{-n+2}\lambda((1-\phi)
{\cal J}_\nu)(\sigma_1\lambda)(1-\phi)(\sigma_2\lambda)$$ $$\times 
\left(e^{i\sigma_2\lambda}\left(b_\nu^+(\sigma_2\lambda)-
c_\nu (\sigma_2\lambda)^{(n-3)/2}\right)+
e^{-i\sigma_2\lambda}\left(b_\nu^-(\sigma_2\lambda)+
c_\nu (\sigma_2\lambda)^{(n-3)/2}\right)\right)$$
$$+c_\nu\sigma_1^{-n+2}\sigma_2^{-(n-1)/2}\lambda^{(n-1)/2}(1-\phi)
(\sigma_1\lambda)(1-\phi)(\sigma_2\lambda)\left(e^{i\sigma_2\lambda}+
e^{-i\sigma_2\lambda}\right)$$ 
$$\times 
\left(e^{i\sigma_1\lambda}\left(b_\nu^+(\sigma_1\lambda)-
c_\nu (\sigma_1\lambda)^{(n-3)/2}\right)-
e^{-i\sigma_1\lambda}\left(b_\nu^-(\sigma_1\lambda)+
c_\nu (\sigma_1\lambda)^{(n-3)/2}\right)\right)$$
 $$=:X(\lambda;\sigma_1,\sigma_2)+Y(\lambda;\sigma_1,\sigma_2)+
Z(\lambda;\sigma_1,\sigma_2),
\eqno{(5.26)}$$
where $X$ denotes the first term in the RHS, $Y$ denotes the second one,
while $Z=Z_1+Z_2$ denotes the remainder ($Z_2$ being the last term). 
In view of (5.25) we have
(for $\lambda^2\in\mbox{supp}\,\psi$)
$$\left|\frac{d^m}{d\lambda^m}Z_1(\lambda;\sigma_1,\sigma_2)
\right|\le C\sum_{k=0}^m\langle\sigma_1\rangle^{k-(n-1)/2}
\langle\sigma_2\rangle^{m-k-(n+1)/2}$$ $$
\le C\langle\sigma_1\rangle^{-(n-1)/2}\langle\sigma_2\rangle^{-(n+1)/2}
\left(\langle\sigma_1\rangle+\langle\sigma_2\rangle\right)^m,
\eqno{(5.27)}$$
for every integer $m\ge 0$.Hence, 
$$\left|\int e^{it\lambda^2}\psi(\lambda^2)Z_1(\lambda;\sigma_1,\sigma_2)
d\lambda\right|\le C|t|^{-m-1/2}\langle\sigma_1
\rangle^{-(n-1)/2}\langle\sigma_2\rangle^{-(n+1)/2}
\left(\langle\sigma_1\rangle+\langle\sigma_2\rangle\right)^m,
\eqno{(5.28)}$$
for all real $m\ge 0$. Take $m=(n-1)/2$ and observe that
$$\langle\sigma_1
\rangle^{-(n-1)/2}\langle\sigma_2\rangle^{-(n+1)/2}
\left(\langle\sigma_1\rangle+\langle\sigma_2\rangle\right)^{(n-1)/2}
\le C\langle\sigma_1\rangle^{-(n+1)/2}+C\langle\sigma_2\rangle^{-(n+1)/2}.$$
Therefore, the integral in the LHS of (5.28) satisfies (5.5). The function
$Z_2$ is treated in precisely the same way. Furthermore, we 
decompose the function
$Y$ as $Y=Y_1+Y_2+Y_3$, with $Y_1$ corresponding to the term 
$\phi(\sigma_1\lambda)$, $Y_2$ corresponding to the term 
$\phi(\sigma_2\lambda)$, and $Y_3$ being the remainder. 
We have (for $\lambda^2\in\mbox{supp}\,\psi$)
$$\left|
\frac{d^m}{d\lambda^m}Y_1(\lambda;\sigma_1,\sigma_2)
\right|
\le C\sigma_1^{-(n-1)/2}\sigma_2^{-(n-1)/2}
\langle\sigma_2\rangle^m,
\eqno{(5.29)}$$
for every integer $m\ge 0$, and hence, 
$$\left|\int e^{it\lambda^2}\psi(\lambda^2)Y_1(\lambda;\sigma_1,\sigma_2)
d\lambda\right|\le C|t|^{-m-1/2}\sigma_1^{-(n-1)/2}
\sigma_2^{-(n-1)/2}\langle\sigma_2\rangle^m,
\eqno{(5.30)}$$
for all real $m\ge 0$. Since $\sigma_1$ is bounded
as long as $\sigma_1\lambda\in\mbox{supp}\,\phi$ and 
$\lambda^2\in\mbox{supp}\,\psi$, we can bound the RHS of (5.30)
(for $m=(n-1)/2$) by 
$$C|t|^{-n/2}\left(\sigma_1^{-(n-1)/2}\sigma_2^{-(n-1)/2}+
\sigma_1^{-(n-1)/2}\right)
\le C|t|^{-n/2}\left(\sigma_1^{-n+1}+\sigma_2^{-n+1}\right).$$
Therefore, the integral in the LHS of (5.30) satisfies (5.5). The terms
corresponding to 
$Y_2$ and $Y_3$ can be treated in precisely the same way.

Furthermore, we write
$$\int e^{it\lambda^2}\psi(\lambda^2)X(\lambda;\sigma_1,\sigma_2)
d\lambda=\int w^{(1)}(\mu,t,\sigma_1,\sigma_2) 
\mu^{(n-3)/4}\psi(\mu)d\mu,\eqno{(5.31)}$$
where the function $w^{(1)}$ is of the form
$$ w^{(1)}(\lambda^2,t,\sigma_1,\sigma_2)
=\sum_\pm\sum_\pm Const(\sigma_1\sigma_2)^{-(n-1)/2}
e^{it\lambda^2+i\lambda(\pm\sigma_1\pm\sigma_2)}\lambda^{(n-3)/2}.
\eqno{(5.32)}$$
Clearly, the function $w^{(1)}$ satisfies (5.4). To prove that
$w^{(1)}$ satisfies (4.10), it suffices to show that 
$$\left|\int_0^\infty
e^{it\lambda^2+i\sigma\lambda}\lambda^{(n-1)/2}\chi_a(\lambda^2)
d\lambda\right|\le C|t|^{-n/2}\langle\sigma\rangle^{(n-1)/2},
\quad\forall t\neq 0,\,\sigma\in{\bf R},\eqno{(5.33)}$$
with a constant $C>0$ independent of $t$ and $\sigma$. 
Consider first the case of $n$ odd, and set 
$m=(n-1)/2$. Integrating by parts $m$ times, we get
$$2(it)^m\int_0^\infty
e^{it\lambda^2+i\sigma\lambda}\lambda^{(n-1)/2}\chi_a(\lambda^2)
d\lambda=(it)^m\int_0^\infty
e^{it\mu+i\sigma\sqrt{\mu}}\mu^{(n-3)/4}\chi_a(\mu)d\mu$$
 $$=\int_0^\infty
e^{it\mu}\frac{d^m}{d\mu^m}\left(e^{i\sigma\sqrt{\mu}}\mu^{(n-3)/4}
\chi_a(\mu)\right)d\mu$$
 $$=\sum_{j=0}^{m-1}\sigma^j\int_0^\infty
e^{it\lambda^2+i\sigma\lambda}\varphi_j(\lambda)d\lambda
+2\int_0^\infty
e^{it\lambda^2+i\sigma\lambda}g_m(\lambda^2,\sigma)\chi_a(\lambda^2)
\lambda d\lambda,\eqno{(5.34)}$$
where $\varphi_j\in C_0^\infty((0,+\infty))$, and 
$$g_m(\mu,\sigma)=e^{-i\sigma\sqrt{\mu}}\frac{d^m}{d\mu^m}
\left(e^{i\sigma\sqrt{\mu}}\mu^{(n-3)/4}\right).$$
By (2.7), each integral in the sum in the RHS of (5.34) is bounded by
$O(|t|^{-1/2})$. To bound the remainder, observe that 
$g_m$ is of the form
$$g_m(\mu,\sigma)=\mu^{-1/2}\sum_{j=0}^m\gamma_j\sigma^j\mu^{-(m-j)/2},$$
where $\gamma_j$ are independent of $\mu$ and $\sigma$. 
Therefore, we have
(for $\lambda^2\in\mbox{supp}\,\chi_a$)
$$\left|\frac{d^j}{d\lambda^j}\left(\lambda g_m(\lambda^2,\sigma)
\right)\right|
\le C\langle\sigma\rangle^m\lambda^{-2j},\quad j=0,1,\eqno{(5.35)}$$
and hence
$$\left|\frac{d}{d\lambda}\left(\lambda g_m(\lambda^2,\sigma)
\chi_a(\lambda^2)\right)\right|
\le C\langle\sigma\rangle^m\langle\lambda\rangle^{-2}.\eqno{(5.36)}$$
We now write the last integral in the RHS of (5.34) as
$$\int_0^\infty k(\lambda,\sigma,t)\frac{d}{d\lambda}\left(\lambda g_m
(\lambda^2,\sigma)\chi_a(\lambda^2)\right)d\lambda,\eqno{(5.37)}$$
where
$$ k(\lambda,\sigma,t)=\int_0^\lambda e^{itx^2+i\sigma x}dx$$
 $$=e^{-i\sigma^2/4t}\int_0^\lambda e^{it(x+\sigma/2t)^2}dx
=|t|^{-1/2}e^{-i\sigma^2/4t}\int_{-\varepsilon\sigma/2|t|^{1/2}}^{
\lambda|t|^{1/2}-\varepsilon\sigma/2|t|^{1/2}}
 e^{i\varepsilon y^2}dy,$$
where $\varepsilon=\mbox{sign}\,t$. Using the well known bound
$$\left|\int_0^a e^{i\varepsilon y^2}dy\right|\le C,\quad \forall a\in
{\bf R},$$
with a constant $C>0$ independent of $a$, we get
$$\left| k(\lambda,\sigma,t)\right|\le C|t|^{-1/2},\eqno{(5.38)}$$
with a constant $C>0$ independent of $\lambda$, $\sigma$ and $t$.
By (5.36) and (5.38), the integral in (5.37) is bounded by
$C\langle\sigma\rangle^m|t|^{-1/2}$, which clearly implies (5.33)
in this case.

Let now $n$ be even and set $m=(n-2)/2$. Then (5.34) still holds and
each integral in the sum in the RHS is bounded by $C_k
\langle\sigma\rangle^k|t|^{-k-1/2}$ for all real $k\ge 0$, and in particular
for $k=1/2$. Therefore, it suffices to show that the last integral
in the RHS of (5.34) is bounded in this case by $C
\langle\sigma\rangle^{m+1/2}|t|^{-1}$. The function $g_m$ in this case 
is of the form
$$g_m(\mu,\sigma)=\mu^{-1/4}\sum_{j=0}^m\gamma'_j\sigma^j\mu^{-(m-j)/2}.$$
Thus, it suffices to show that 
$$\left|\int_0^\infty
e^{it\lambda^2+i\sigma\lambda}\lambda^{1/2-j}\chi_a(\lambda^2)
d\lambda\right|\le C|t|^{-1}\langle\sigma\rangle^{1/2+j},
\quad 0\le j\le m.\eqno{(5.39)}$$
When $j\ge 1$, (5.39) follows easily by integrating once by parts. To prove
(5.39) for $j=0$, we proceed as follows (if $\sigma/t<0$)
$$\int_0^\infty
e^{it\lambda^2+i\sigma\lambda}\lambda^{1/2}\chi_a(\lambda^2)
d\lambda=e^{-i\sigma^2/4t}\int_0^\infty
e^{it(\lambda+\sigma/2t)^2}\lambda^{1/2}\chi_a(\lambda^2)
d\lambda$$
 $$=(-\sigma/2t)^{1/2}e^{-i\sigma^2/4t}\int_0^\infty
e^{it(\lambda+\sigma/2t)^2}\chi_a(\lambda^2)
d\lambda$$ $$+
e^{-i\sigma^2/4t}\int_0^\infty
e^{it(\lambda+\sigma/2t)^2}\frac{\chi_a(\lambda^2)}{\lambda^{1/2}+
(-\sigma/2t)^{1/2}}
d(\lambda+\sigma/2t)^2$$
$$=(-\sigma/2t)^{1/2}e^{-i\sigma^2/4t}\int_0^\infty
e^{it(\lambda+\sigma/2t)^2}\chi_a(\lambda^2)
d\lambda$$ $$-(it)^{-1}e^{-i\sigma^2/4t}\int_0^\infty
e^{it(\lambda+\sigma/2t)^2}\frac{d}{d\lambda}\left(
\frac{\chi_a(\lambda^2)}{\lambda^{1/2}+
(-\sigma/2t)^{1/2}}\right)d\lambda. \eqno{(5.40)}$$
We bound the integral in the first term in the RHS of (5.40) by
$O(|t|^{-1/2})$ in the same way as the integral (5.37) above, so the first
term itself is bounded by $C\langle\sigma\rangle^{1/2}|t|^{-1}$.
The second term is bounded by $O(|t|^{-1})$ because of the bound
$$\left|\frac{d}{d\lambda}\left(
\frac{\chi_a(\lambda^2)}{\lambda^{1/2}+
(-\sigma/2t)^{1/2}}\right)\right|\le C\langle\lambda\rangle^{-3/2},$$
with a constant $C>0$ independent of $\lambda$, $\sigma$ and $t$. 
When $\sigma/t\ge 0$, we write 
$$\int_0^\infty
e^{it\lambda^2+i\sigma\lambda}\lambda^{1/2}\chi_a(\lambda^2)
d\lambda$$ $$=-(it)^{-1}e^{-i\sigma^2/4t}\int_0^\infty
e^{it(\lambda+\sigma/2t)^2}\frac{d}{d\lambda}\left(
\frac{\lambda^{1/2}\chi_a(\lambda^2)}{\lambda+\sigma/2t}
\right)d\lambda, \eqno{(5.41)}$$
so the integral in the LHS of (5.41) 
is bounded by  $O(|t|^{-1})$ because of the bound
$$\left|\frac{d}{d\lambda}\left(
\frac{\lambda^{1/2}\chi_a(\lambda^2)}{\lambda+\sigma/2t}
\right)\right|\le C\langle\lambda\rangle^{-3/2},$$
with a constant $C>0$ independent of $\lambda$, $\sigma$ and $t$. 
This completes the proof of (5.33). Since the function 
$U_1^{(2)}$ can be treated in precisely the same way as $U_1^{(1)}$,
the proof of the proposition is completed.

G. Vodev, Universit\'e de Nantes,
 D\'epartement de Math\'ematiques, UMR 6629 du CNRS,
 2, rue de la Houssini\`ere, BP 92208, 44332 Nantes Cedex 03, France

e-mail: georgi.vodev@math.univ-nantes.fr


\begin{thebibliography} 
\frenchspacing \baselineskip=12 pt plus 1pt minus 1pt 

\bibitem{kn:G} {\sc M. Goldberg}, {\em Dispersive bounds for the
three-dimensional Schr\"odinger equation with almost critical potentials},
GAFA, to appear.

\bibitem{kn:GS} {\sc M. Goldberg and W. Schlag}, {\em Dispersive estimates
 for Schr\"odinger operators in dimensions one and three}, 
Commun. Math. Phys. {\bf 251} (2004), 157-178.

\bibitem{kn:GV} {\sc M. Goldberg and M. Visan}, {\em A counterexample 
to dispersive estimates for Schr\"odinger operators 
in higher dimensions}, preprint. 

\bibitem{kn:JN} {\sc A. Jensen and S. Nakamura}, {\em $L^p$-mapping
properties of functions of  Schr\"odinger operators and their
applications to scattering theory}, J. Math. Soc. Japan {\bf 47} (1995),
253-273.

\bibitem{kn:JSS} {\sc J.-L. Journ\'e, A. Sofer and C. Sogge}, {\em 
 Decay estimates for Schr\"odinger operators}, Commun. Pure Appl. Math.
{\bf 44} (1991), 573-604.

\bibitem{kn:RS} {\sc I. Rodnianski and W. Schlag}, {\em Time decay for
solutions of Schr\"odinger equations with rough and time
-dependent potentials}, Invent. Math. {\bf 155} (2004), 451-513.

\bibitem{kn:S} {\sc W. Schlag}, {\em Dispersive estimates
 for Schr\"odinger operators in two dimensions}, 
Commun. Math. Phys. {\bf 257} (2005), 87-117.

\bibitem{kn:V1} {\sc G. Vodev},  
{\em Dispersive estimates of solutions to the Schr\"odinger equation},
Ann. Henri Poincar\'e {\bf 6} (2005), 1179-1196.

\bibitem{kn:V2} {\sc G. Vodev},  
{\em Dispersive estimates of solutions to the wave equation with a potential 
in dimensions $n\ge 4$}, Commun. Partial Diff. Equations, to appear.

\bibitem{kn:Y1} {\sc K. Yajima}, {\em The $W^{k,p}$-continuity of wave
operators for Schr\"odinger operators}, J. Math. Soc. Japan {\bf 47} (1995),
551-581.

\bibitem{kn:Y2} {\sc K. Yajima}, {\em Dispersive estimates
 for Schr\"odinger equations with threshold resonance and
eigenvalue}, Commun. Math. Phys. {\bf 259} (2005), 475-509.

\end{thebibliography}
\end{document}